\documentclass[11pt,a4paper]{amsart}
\usepackage{upref}
\usepackage{bm}
\usepackage{mathrsfs}
\usepackage{marginnote}

\usepackage{lmodern} 
\usepackage{amsmath}
\usepackage{amssymb} 
\usepackage{stmaryrd} 
\usepackage{mathtools} 
\usepackage{mathrsfs}  
\usepackage{graphicx} 
\usepackage{xcolor}   
\usepackage[english]{babel} 
\usepackage{caption} 
\usepackage{subcaption} 
\usepackage{hyperref} 
\usepackage{thmtools} 
\usepackage{microtype} 
\usepackage[shortlabels]{enumitem} 
\usepackage{xargs} 
\usepackage[colorinlistoftodos,prependcaption,textsize=scriptsize]{todonotes} 
\usepackage{csquotes} 
\usepackage{dsfont} 
\usepackage[title]{appendix}
\usepackage{subfiles}
\usepackage{import}
\usepackage{url}

\addto\extrasenglish{%
}

\numberwithin{equation}{section}

\newtheorem{theorem}{Theorem}[section]
\newtheorem{proposition}[theorem]{Proposition}
\newtheorem{lemma}[theorem]{Lemma}
\newtheorem{corollary}[theorem]{Corollary}
\theoremstyle{definition}
\newtheorem{definition}[theorem]{Definition}

\newtheorem{example}[theorem]{Example}
\theoremstyle{remark}
\newtheorem{remark}[theorem]{Remark}
\newtheorem{notation}[theorem]{Notation}

\newcommand\R{{\mathbb{R}}}
\newcommand\N{\mathbb{N}}
\DeclareMathOperator{\supp}{supp}

\newcommandx{\lucas}[2][1=]{\todo[linecolor=red,backgroundcolor=red!25,bordercolor=red,#1, author = Lucas]{#2}} 
\newcommandx{\david}[2][1=]{\todo[linecolor=blue,backgroundcolor=blue!25,bordercolor=blue,#1, author = David]{#2}} 

\mathtoolsset{showonlyrefs} 

\title[Besov Reconstruction]{Besov Reconstruction}

\author{Lucas Broux}
\address[Lucas Broux]{Sorbonne Universit\'e, Laboratoire de
	Probabilit\'es Statistique et Mod\'elisation, 4 Pl. Jussieu, 75005 Paris,France}
\email{\href{mailto:lucas.broux@upmc.fr }{\nolinkurl{lucas.broux@upmc.fr }}}
\author{David Lee}
 \address[David Lee]{Sorbonne Universit\'e, Laboratoire de
   Probabilit\'es Statistique et Mod\'elisation, 4 Pl. Jussieu, 75005 Paris,France}
\email{\href{mailto:david.lee@upmc.fr}{\nolinkurl{david.lee@upmc.fr}}}

\makeatletter
\@namedef{subjclassname@2020}{%
  \textup{2020} Mathematics Subject Classification}
\makeatother


\begin{document}

\begin{abstract}
The reconstruction theorem tackles the problem of building a global distribution, on $\mathbb{R}^d$ or on a manifold, for a given family of sufficiently coherent local approximations. 
This theorem is a critical tool within Hairer's theory of Regularity Structures.
In this paper, we establish a reconstruction theorem in the Besov setting, extending recent results of Caravenna and Zambotti.
A Besov reconstruction theorem was first formulated by Hairer and Labb\'e in the context of regularity structures, exploiting nontrivial results from wavelet analysis. 
Our calculations follow the more elementary approach of coherent germs due to Caravenna and Zambotti.
With this formulation our results are both stated and proved with tools from the theory of distributions without the need of the theory of Regularity Structures.
As an application, we present an alternative proof of a (Besov) Young multiplication theorem which does not require the use of para-differential calculus.
\end{abstract}

\maketitle

\tableofcontents


\section{Motivation and Background}


Multiplying two distributions in a general setting is a notoriously difficult problem in many situations in PDEs and mathematical analysis.
Recently this problem has motivated an intensive activity in stochastic analysis. Inspired by the theory of Rough Paths \cite{MR1654527}, two approaches to singular SPDEs have been developed in the last decade: regularity structures
\cite{MR3274562, MR3935036, chandra2018analytic, MR4210726} and paracontrolled distributions \cite{MR3406823}. These new theories allow to give a meaning and also a well-posedness
result for PDEs with stochastic forcing terms which would be ill-defined using classical tools; indeed, due to the wild oscillations of the noise, the solutions are expected to be distributions; if the equations contain polynomial or analytical non-linearities, then such terms are ill-defined
and require new ideas.

It is not our aim to enter into the details of regularity structures. On the contrary, we want to present one of its main results, the \emph{reconstruction theorem} in Besov spaces, in a 
more elementary way. 

One of the main ideas of the theory is to lift the equation to a space of \enquote{local approximations} of the solutions, rather than to work directly in the space of Schwartz distributions. The reconstruction theorem allows to retrieve a genuine distribution from such a family of local approximations.
Note that although the theory permits deep results in stochastic analysis, this theorem is purely deterministic.

More precisely, we consider the following \enquote{reconstruction problem}:
\begin{center}\emph{
	Given the data, for all $x \in \mathbb{R}^d$, of a distribution $F_x \in \mathcal{D}^{\prime} ( \mathbb{R}^d )$, is there a distribution $f$ on $\R^d$ which is well approximated by $F_{x}$ around each point $x\in \R^d$?}
\end{center} 

Of course, if for all $x \in \mathbb{R}^d$, $F_x$ is an actual continuous function, then the function $f \colon x \mapsto F_x(x)$ is a natural answer to this question. 

For instance, if $F_x(\cdot)$ is the Taylor polynomial of order $r \in \mathbb{N}$, at the point $x \in \mathbb{R}^d$, of some smooth function $f \colon \mathbb{R}^d \to \mathbb{R}$, then it indeed holds that $f \left( x \right) = F_x(x)$; furthermore $f$ is well approximated by $F_x$ around $x$, with the following estimate: $|f(y)-F_{x}(y)|\lesssim |y-x|^{r + 1}$ for $y$ close to $x$.
Note that in this situation one also has $| F_z \left( w \right) - F_y \left( w \right) | \lesssim \sum_{| l | < r} | w - y |^{|l|} | y - z |^{r - | l |}$, which asserts that the family $(F_x)_{x \in \mathbb{R}^d}$ is sufficiently \enquote{coherent}.

The situation becomes more subtle when the objects at play are actual distributions, which is the case in the previously mentioned context of stochastic PDEs, where distributional terms arise from the white noise governing the equation.

The reconstruction theorem \cite[Theorem~3.10]{MR3274562}, which solves the above problem in the context of H\"older spaces, was originally stated in the formalism of regularity structures, and proved using wavelet analysis. 
Alternative proofs, within the context of regularity structures, were later established in \cite{MR3406823, MR3925533, MR4174393}.

More recently, the theorem was revisited in \cite[Theorem~5.1]{caravenna2020hairers}, where it was stated and proven in an elementary 
and more general setting.
In particular, \cite{caravenna2020hairers} exhibits a sufficient condition on $\left( F_x \right)_{x \in \mathbb{R}^d}$, dubbed \emph{coherence}, under which such a reconstruction $f$ exists.
The reconstruction $f$ is then built from a custom-made dyadic decomposition inspired by mollification, which replaces the wavelet decomposition of the original theorem.
See also \cite{rinaldi2021reconstruction} for a version of this result over smooth manifolds (still in the H\"older setting).

Unfortunately, it turns out that working in H\"older spaces may not be enough for many purposes, and it is desirable to consider the more general Besov spaces $\mathcal{B}_{p, q}^{\alpha} (\mathbb{R}^d)$ with $p,q\in [1,\infty]$, as they allow for finer analysis of distributions.
For instance, it is natural to consider the Dirac mass $\delta_0$ on $\mathbb{R}^d$ as an initial condition for some stochastic PDEs.
In such cases, one may benefit from the fact that not only does $ \delta_0 \in \mathcal{C}^{- d}$, it also holds that $\delta_0 \in \mathcal{B}_{p, \infty}^{- d + d / p}$ for any $p \in \left[ 1, + \infty \right]$, thus allowing to work with improved regularity (at the expense of integrability); see \cite{MR3779690} for an application of this idea.

Another useful property of Besov spaces is that when $p=q=2$, they match with (fractional) Sobolev spaces, which are the natural framework in many situations.
For instance, this allows to construct random differential operators and study their spectral properties \cite{MR3997634}, or to apply Malliavin calculus to solutions of stochastic PDEs \cite{MR4058990}. 

This motivates the need for a reconstruction theorem in the more general context of Besov spaces, which, as it turns out, has already been established in the formalism of regularity structures in \cite[Theorem 3.1]{MR3684891}, using once again wavelets in its proof, and later in \cite{singh2018elementary}.
See also \cite{liu2021sobolev}, where a similar reconstruction result is proposed and applied to the problem of lifting Sobolev paths to Sobolev rough paths.
Note that a Besov sewing lemma -- an analogous result in the context of Rough Paths theory -- has recently been established in \cite{friz2021besov} and applied to rough differential equations.

In this article, we provide a version of this Besov reconstruction theorem, in the spirit of \cite{caravenna2020hairers}.
In particular, not only does our result generalise both \cite[Theorem 3.1]{MR3684891} and \cite[Theorem~5.1]{caravenna2020hairers} (in the case of global exponents), it is independent of the theory of regularity structures and can be formulated in the language of distributions.
Note that we will not talk about regularity structures in this paper, except in comparing our results to \cite{MR3684891}.

We present our main results, \autoref{corollary:reconstruction_with_coherence_homogeneity} and \autoref{Prop:Existence_gamma_gq_0}, as sufficient conditions on the family $\left( F_x \right)_{x \in \mathbb{R}^d}$ for a reconstruction to exist in a prescribed Besov sense. 
Interestingly, in \autoref{corollary:reconstruction_with_coherence_homogeneity}, we exhibit simple notions of \enquote{coherence} and \enquote{homogeneity}, which generalise the results of \cite[Theorem 3.1]{MR3684891} and \cite[Theorem~5.1]{caravenna2020hairers}.
It actually turns out that the conditions of \autoref{corollary:reconstruction_with_coherence_homogeneity} can be refined, and we propose a more general reconstruction result as \autoref{Prop:Existence_gamma_gq_0}, which we discuss in a later section because its statement is more technical.
Remarkably, \autoref{Prop:Existence_gamma_gq_0} allows us to tackle the classical problem of constructing a product between suitable Besov spaces, as stated in \autoref{prop:Young_multiplication_besov}.
Contrary to usual approaches, \autoref{prop:Young_multiplication_besov} does not require paraproducts, and to the best of our understanding our conditions on the parameters of the spaces are quite optimal, see \autoref{section:application_young_multiplication} for a short review of the literature.
See also \cite[Section~14]{caravenna2020hairers} for a similar application in the H\"older case.

\subsection*{Outline}
This paper is organised as follows. 
In \autoref{section:notations}, we set the main notations that will be used in the remainder of the paper.

In \autoref{section:main_results}, we discuss the problem of reconstruction that we consider, state an important reconstruction result as \autoref{corollary:reconstruction_with_coherence_homogeneity}, and compare it to the results of \cite{caravenna2020hairers,MR3684891}.
We also discuss the problem of building a Young multiplication in Besov spaces, that is, a continuous bilinear map extending the usual pointwise product between smooth functions. We construct such a product in \autoref{prop:Young_multiplication_besov}.

\autoref{section:proof_of_main_result} is devoted to the proof of a general reconstruction result, \autoref{Prop:Existence_gamma_gq_0}, following an approach similar to \cite{caravenna2020hairers}.

In \autoref{section:proof_of_corollary}, we prove \autoref{corollary:reconstruction_with_coherence_homogeneity} and \autoref{thm:properties_of_the_reconstruction_map} as a corollary of  \autoref{Prop:Existence_gamma_gq_0}.

Finally, in \autoref{section:application_young_multiplication}, we prove \autoref{prop:Young_multiplication_besov} as a corollary of  \autoref{Prop:Existence_gamma_gq_0}.


\section{Notations}
\label{section:notations}

In this paper, we consider an integer $d \in \mathbb{N}$ and work in the space $\mathbb{R}^d$, equipped with its canonical Euclidean norm $\left| \cdot \right|$.
If $x \in \mathbb{R}^d$ and $r \in \mathbb{R}_+$, we will denote $B \left( x, r \right) \coloneqq \lbrace z \in \mathbb{R}^d, \left| z - x \right| \leq r \rbrace$ to be the closed ball of center $x$ and radius $r$.
When $K \subset \mathbb{R}^d$, and $R > 0$ we will denote $\bar{K}_R \coloneqq K + B (0, R)$ its $R$-enlargement.

We shall work extensively in the following Lebesgue spaces, where the variables and domain of integration will usually be clear from context:
	\begin{enumerate}
		\item The variable $x \in \mathbb{R}^d$ will usually correspond to a space variable and we denote:
			\begin{equation}
				L^p = L^p \left( x \right) \coloneqq L^p \left( \mathbb{R}^d,d x \right) .
			\end{equation}
		\item The variable $h \in \mathbb{R}^d$ will usually correspond to a space variable and we denote:
			\begin{equation}
				L_{h}^q = L_h^q \left( h \right) \coloneqq L^q \left(B(0, 1), \frac{d h}{\left| h \right|^d} \right) .
			\end{equation}
			
			We might also integrate on a domain $K$ rather than $B \left( 0, 1 \right)$. In this case, we will denote:
				\begin{equation}
					L_h^q \left( h \in K \right) \coloneqq L^q \left(K, \frac{d h}{\left| h \right|^d} \right) .
				\end{equation}
			
		\item The variable $\lambda \in \left( 0, 1 \right]$ will usually correspond to a scaling variable and we denote:
			\begin{equation}
				L_{\lambda}^q = L_{\lambda}^q \left( \lambda \right) \coloneqq L^{q} \left( \left( 0, 1 \right] ,\frac{d \lambda}{\lambda} \right)
			\end{equation}
		
		\item The variable $n \in \mathbb{N}$ will usually correspond to a scaling variable and we denote:
			\begin{equation}
				\ell^q = \ell^q \left( n \right) \coloneqq \ell^q \left( n \in \mathbb{N} \right) .
			\end{equation}
			We might also sum on $n \geq n_0$ rather than $\mathbb{N}$. For $n_0 \in \mathbb{Z}$, we will denote $\ell^q \left( n \geq n_0 \right)$ the corresponding space.
	\end{enumerate}
In the case where $p$ or $q$ are equal to $\infty$ then the associated norm  corresponds to the usual supremum norm.

The space of test-functions is denoted $\mathcal{D} = \mathcal{D} ( \mathbb{R}^d )$ and is defined as the space of $\mathcal{C}^{\infty} ( \mathbb{R}^d )$ functions with compact support.
More generally, if $K$ is a compact set of $\mathbb{R}^d$, $\mathcal{D} ( K )$ will denote the space of test-functions supported in $K$.

If $\varphi \colon \mathbb{R}^d \to \mathbb{R}$ is a sufficiently differentiable function, we will denote its partial derivatives by, for a multi-index $k = \left( k_1, \cdots , k_d \right) \in \mathbb{N}^d$: $\partial^k \varphi \coloneqq \partial_1^{k_1} \cdots \partial_d^{k_d} \varphi$.

Let $r \in \mathbb{N}$, the $\mathcal{C}^r$ norm of a sufficiently differentiable function $\varphi$ is defined by:
	\begin{equation}
		\left\| \varphi \right\|_{\mathcal{C}^r} \coloneqq \max\limits_{\left| k \right| \leq r} \left\| \partial^k \varphi \right\|_{\infty} .
	\end{equation}

Let $\varphi \in \mathcal{D} ( \mathbb{R}^d )$, $x \in \mathbb{R}^d$, $\lambda > 0$, we denote $\varphi_x^{\lambda}$ the scaled and recentered version of $\varphi$, defined as follows: $\varphi_x^{\lambda} \left( \cdot \right) \coloneqq \lambda^{- d} \varphi \left( \lambda^{-1} \left( \cdot - x \right) \right)$.

We will denote the multinomial $x^k = \prod_i x_i^{k_i}$. Often we will use the notation $|k|:=k_1+...+k_d$. 

For $r, s \in \mathbb{N}$ and $K$ a compact set of $\mathbb{R}^d$ we define:

	\begin{equation}
		\begin{dcases}
			\mathscr{B}^r \left( K \right) & \coloneqq \left\lbrace \varphi \in \mathcal{D} \left( K \right), \left\| \varphi \right\|_{\mathcal{C}^r} \leq 1 \right\rbrace , \\
			\mathscr{B}_s^r \left( K \right) & \coloneqq \left\lbrace \varphi \in \mathcal{D} \left( K \right), \left\| \varphi \right\|_{\mathcal{C}^r} \leq 1 \text{ and } \int x^k \varphi \left( x \right) d x = 0 , 0 \leq \left| k \right| \leq s \right\rbrace .
		\end{dcases}
	\end{equation}

Most of the time, $K$ will be $B \left( 0, 1 \right)$, so for simplicity of notation, we will denote:
	\begin{equation}
		\begin{dcases}
			\mathscr{B}^r & \coloneqq \mathscr{B}^r \left( B \left( 0, 1 \right) \right) , \\
			\mathscr{B}_s^r & \coloneqq \mathscr{B}_s^r \left( B \left( 0, 1 \right) \right) .
		\end{dcases}
	\end{equation}

Recall that a (Schwartz) distribution is a linear functional $f \colon \mathcal{D} ( \mathbb{R}^d ) \to \mathbb{R}$ such that for all compact $K \subset \mathbb{R}^d$, there exists $r = r_K \in \mathbb{N}$ and $C = C_K < + \infty$ such that for all $\varphi \in \mathcal{D} ( \mathbb{R}^d )$ supported in $K$, $\left| f \left( \varphi \right) \right| \leq C \left\| \varphi \right\|_{\mathcal{C}^r}$.
When $r$ does not depend on $K$, we say that $f$ is a distribution of order $r$.
We denote $\mathcal{D}^{\prime} ( \mathbb{R}^d )$ the space of distributions.
When $f \colon \mathcal{D} ( \mathbb{R}^d ) \to \mathbb{R}$ is a linear functional satisfying the above condition, for a given compact set $K$, we say that $f$ is a distribution on $K$ and we note the corresponding space $\mathcal{D}^{\prime} ( K )$.

We will denote $\mathcal{B}_{p, q}^{\alpha}  = \mathcal{B}_{p, q}^{\alpha} ( \mathbb{R}^d )$ the (nonhomogeneous) Besov spaces of exponent $\alpha \in \mathbb{R}$ and integrability parameters $p, q \in \left[ 1, + \infty \right]$, on the whole space $\mathbb{R}^d$. 
We will denote $\mathcal{B}_{p, q, \mathrm{loc}}^{\alpha}$ the corresponding local Besov space.
See \autoref{appendix_besov_spaces} for the definition and properties of Besov spaces.

The reader might be surprised that we work in $\mathcal{D}^{\prime}$ rather than in the space of tempered distributions $\mathcal{S}^{\prime}$, which is more natural in the context of Besov spaces.
However, this is not really problematic, see \autoref{remark:Besov_spaces_tempered_distributions} below.

\section{Main Results}
\label{section:main_results}

\subsection{The problem of reconstruction}

To begin, let us properly define the notion of a \emph{germ}. 

\begin{definition}[Germ]
A \emph{germ} is a family of distributions $\left( F_x \right)_{x \in \mathbb{R}^d}$, i.e.\ for all $x \in \mathbb{R}^d$, $F_x \in \mathcal{D}^{\prime} ( \mathbb{R}^d )$.
We also assume that for all test-functions $\varphi \in \mathcal{D} ( \mathbb{R}^d )$, the map $x \mapsto F_x \left( \varphi \right)$ is measurable.

For simplicity, we will also denote $(F_{x})_{x\in \R^d}$ as $F$.
\end{definition}

We think of a germ as a family of local approximations for a global distribution $\mathcal{R} ( F )$ that is to be reconstructed in a suitable Besov sense.
For this purpose, we shall consider in our main result the following scaling functions.
Fix $\epsilon > 0$ arbitrary.
For $\gamma \in \mathbb{R}$, $q \in \left[ 0, + \infty \right]$, and $\lambda \in \left( 0, 1 \right]$, set:
	\begin{equation}
		k \left( \lambda \right) \coloneqq k_{\gamma, q, \epsilon} \left( \lambda \right) \coloneqq
		\begin{dcases}
			\lambda^{\gamma} & \text{if } \gamma \neq 0 , \\
			1 + \left| \log \left( \lambda \right) \right| & \text{if } \gamma = 0, q=+\infty , \\
			1 + \left| \log \left( \lambda \right) \right|^{1 + \epsilon} & \text{if } \gamma = 0, q<+\infty .
		\end{dcases}
		\label{eq:definition_scaling_function}
	\end{equation}

One important result that we prove in this paper (but not our most general, see \autoref{remark:this_result_is_not_the_most_general} below) is the following.

\begin{theorem}[Besov reconstruction]
\label{corollary:reconstruction_with_coherence_homogeneity}
Let $F$ be a germ, $p, q \in \left[ 1, + \infty \right]$, and $\alpha, \beta, \gamma \in \mathbb{R}$ be such that $\alpha \leq \gamma$.
Assume that there exists a test-function $\varphi$ such that $\int \varphi \neq 0$ and that for all $K \subset \mathbb{R}^d$, the following \enquote{homogeneity} property is satisfied:
	\begin{align}
			\left\| F \right\|_{p, \beta, K, \varphi}^{\mathrm{hom}} & \coloneqq \left\| \left\| \frac{F_x \left( \varphi_x^{2^{-n}} \right)}{2^{- n \beta}} \right\|_{L^{p} \left( K, d x \right)} \right\|_{\ell^{\infty} \left( n \in \mathbb{N} \right)} < + \infty .
			\label{eq:homogeneity_condition}
	\end{align}
Assume also that the following \enquote{coherence} property is satisfied:
\begin{align}
	\left\| F \right\|_{p, q, \alpha, \gamma, K, \varphi}^{\mathrm{coh}} & \coloneqq  \left\| \left\| \left\| \frac{\left( F_{x + h} - F_x \right) \left( \varphi_x^{2^{-n}} \right)}{2^{- n \alpha} \left( 2^{-n} + \left| h \right| \right)^{\gamma - \alpha}} \right\|_{L^{p} \left( K, d x \right)} \right\|_{\ell^{\infty} \left( n \in \mathbb{N} \right)} \right\|_{L^q \left( B \left( 0, 2 \right), \frac{dh}{| h |^d} \right)} \\
	& < + \infty . \label{eq:coherence_condition_1}
\end{align}

Then there exists $\mathcal{R} \left( F \right) \in \mathcal{D}^{\prime} ( \mathbb{R}^d )$ satisfying the following reconstruction bound for any integer $r > \max ( -\alpha, - \beta )$  and any $K \subset \mathbb{R}^d$ (recall that $k$ is defined in \eqref{eq:definition_scaling_function}):
		\begin{equation} \label{eq:reconstruction_bound_main_result}
			\left\| \left\| \sup\limits_{\psi \in \mathscr{B}^{r}} \left| \frac{\left( \mathcal{R} \left( F \right) - F_x \right) \left( \psi_x^{\lambda} \right)}{k \left( \lambda \right)} \right| \right\|_{L^{p} \left( K, d x \right)} \right\|_{L^q ( B \left( 0, 1 \right), \frac{d \lambda}{\lambda} )} < + \infty .
		\end{equation}

In fact, the quantity on the left-hand side of \eqref{eq:reconstruction_bound_main_result} can be bounded by a constant times $\left\| F \right\|_{p, q, \alpha, \gamma, \bar{K}_2, \varphi}^{\mathrm{coh}}$.
Furthermore, such an $\mathcal{R} \left( F \right)$ is unique when $\gamma > 0$ but not when $\gamma \leq 0$.
\end{theorem}

Let us propose a few remarks before giving more precise results on the reconstruction map.

\begin{remark}\label{remark:this_result_is_not_the_most_general}
\autoref{corollary:reconstruction_with_coherence_homogeneity} is interesting because, as we will establish in \autoref{ex:Caravenna_Zambotti} and \autoref{ex:Hairer_Labbe}, it is a generalisation of \cite[Theorem 3.1]{MR3684891} and \cite[Theorem~5.1]{caravenna2020hairers} (in the case of global exponents).
However, note that it is not clear whether the condition \eqref{eq:coherence_condition_1} is canonical, in particular regarding the order of integration in the Lebesgue norms.
This suggests that \eqref{eq:coherence_condition_1} is not the \enquote{optimal} condition.
Indeed, we actually establish a more general result of reconstruction, which we state in \autoref{Prop:Existence_gamma_gq_0}.
However, the conditions of \autoref{Prop:Existence_gamma_gq_0} require heavier notations, which is why we postpone its statement and proof to \autoref{section:proof_of_main_result}.
In the remainder of this paper, we shall first prove \autoref{Prop:Existence_gamma_gq_0}, in \autoref{section:proof_of_main_result}.
Then, in \autoref{section:proof_of_corollary}, we shall establish \autoref{corollary:reconstruction_with_coherence_homogeneity} as a corollary of \autoref{Prop:Existence_gamma_gq_0}.
\end{remark}

\begin{remark}
The case $\gamma = 0$ appears as a critical case. Recall that this is the case also in \cite[cf.~Theorem~5.1]{caravenna2020hairers}.
\end{remark}

\begin{remark}
Note that in condition \eqref{eq:coherence_condition_1}, we require integration over $h \in B \left( 0, 2 \right)$, while in the reconstruction bound \eqref{eq:reconstruction_bound_main_result} one integrates over $\lambda \in B \left( 0, 1 \right)$.
It would be more natural to impose the constraint of integrability over $h \in B \left( 0, 1 \right)$, as one would expect the estimates to propagate from $B \left( 0, 1 \right)$ to $B \left( 0, 2 \right)$, similarly to \cite{caravenna2020hairers}.
However, in our context, the asymmetry between the roles of the variables $x$ and $h$, and the fact that the variables $x, h, n$ are \enquote{linked} by the integration, prevent the same argument as in \cite{caravenna2020hairers} to be applied. Of course, in practical situations it is usually equivalent to check the conditions for $B(0,1)$ or for $B(0,2)$.
\end{remark}

\begin{remark}
Note also that contrary to \cite{caravenna2020hairers}, here we do not consider the converse problem of whether the existence of a reconstruction implies any coherence condition such as \eqref{eq:coherence_condition_1}.
\end{remark}

We shall also show that the reconstruction map $\mathcal{R}$ admits the following properties (and see \autoref{Prop:Existence_gamma_gq_0} for a more general version).

\begin{theorem}[Properties of the reconstruction map]\label{thm:properties_of_the_reconstruction_map}
In the context of \autoref{corollary:reconstruction_with_coherence_homogeneity}:
	\begin{enumerate}
		\item (Global version) if \eqref{eq:homogeneity_condition}, \eqref{eq:coherence_condition_1} are satisfied for $K = \mathbb{R}^d$, then:
		\begin{enumerate}
			\item If $\beta \wedge \gamma > 0$, then $\mathcal{R} ( F ) = 0$.
			\item If $\beta \wedge \gamma = 0$, then for all $\kappa > 0$, $\mathcal{R} \left( F \right) \in \mathcal{B}_{p, 1}^{- \kappa}$.
			\item If $\beta \wedge \gamma < 0$, then $\mathcal{R} \left( F \right) \in \mathcal{B}_{p, \infty}^{\beta \wedge \gamma}$.
		\end{enumerate}
	
	\item (Local version) if \eqref{eq:homogeneity_condition}, \eqref{eq:coherence_condition_1} are satisfied for all $K \subset \mathbb{R}^d$, then:
		\begin{enumerate}
			\item If $\beta \wedge \gamma > 0$, then $\mathcal{R} ( F ) = 0$.
			\item If $\beta \wedge \gamma = 0$, then for all $\kappa > 0$, $\mathcal{R} \left( F \right) \in \mathcal{B}_{p, 1, \mathrm{loc}}^{- \kappa}$.
			\item If $\beta \wedge \gamma < 0$, then $\mathcal{R} \left( F \right) \in \mathcal{B}_{p, \infty, \mathrm{loc}}^{\beta \wedge \gamma}$.
		\end{enumerate}
	\end{enumerate}

The reconstruction map is continuous in the following sense: let $\mathcal{B}$ denote the Besov space $\mathcal{B}_{p, 1, K}^{- \kappa}$ if $\beta \wedge \gamma = 0$ and $\mathcal{B}_{p, \infty, K}^{\beta \wedge \gamma}$ if $\beta \wedge \gamma < 0$, then:
	\begin{equation}
		\left\| \mathcal{R} \left( F \right) \right\|_{\mathcal{B}} \lesssim \left\| F \right\|_{p, \beta, \bar{K}_2, \varphi}^{\mathrm{hom}} + \left\| F \right\|_{p, q, \alpha, \gamma, \bar{K}_4, \varphi}^{\mathrm{coh}} .
	\end{equation}
\end{theorem}

\subsection{A comparison with the literature}

Now let us compare \autoref{corollary:reconstruction_with_coherence_homogeneity} and \autoref{thm:properties_of_the_reconstruction_map} above with the existing literature.

\begin{example}[Caravenna-Zambotti]\label{ex:Caravenna_Zambotti}
Taking $p = q = + \infty$ in the previous theorem (in its local version) retrieves \cite[Theorem~5.1]{caravenna2020hairers} in the situation where the coherence and homogeneity exponents $\alpha$ and $\beta$ do not depend on the compact $K$.
Recall also that in the context of \cite{caravenna2020hairers}, the property of (local) homogeneity is implied by the property of (local) coherence, cf.\ \cite[Lemma~4.12]{caravenna2020hairers}.
However, it is not clear whether this generalises to the case of global exponents i.e.\ when $\alpha, \beta, \gamma$ do not depend on the choice of $K$, which is why we assume both homogeneity and coherence.
\end{example}

\begin{example}[Hairer-Labb\'e]
\label{ex:Hairer_Labbe}
Let us shortly discuss how \autoref{corollary:reconstruction_with_coherence_homogeneity} (in its global version) generalises \cite[Theorem~3.1]{MR3684891}, in the case of the canonical scaling $\mathfrak{s} = \left( 1, \dots, 1 \right)$.
Here we assume that the reader is familiar with the framework and notations of \cite{MR3684891}.
Let $(\mathscr{A},\mathscr{T},\mathscr{G})$ be a regularity structure over $\mathbb{R}^d$, endowed with a model $(\Pi,\Gamma)$.
Let $\gamma > 0$, $p, q \in \left[ 1, + \infty \right]$ and let $f\in \mathcal{D}^{\gamma}_{p,q}$ be a Besov modelled distribution. We define a germ $F$ by setting for $x \in \mathbb{R}^d$, $F_x \coloneqq \Pi_x \left( f \left( x \right) \right)$.
Then we claim that $F$ satisfies the homogeneity property \eqref{eq:homogeneity_condition} and the coherence property \eqref{eq:coherence_condition_1}, so that $\mathcal{R} \left( F \right)$ coincides with $\mathcal{R}^{\rm{HL}} \left( f \right)$,
the Hairer-Labbé reconstruction of $f$. 
Let us only discuss the coherence, as the homogeneity is obtained with a similar argument.	
One can observe that for any test-function $\varphi \in \mathcal{D} ( \mathbb{R}^d )$,
	\begin{align*}
	(F_{x+h}-F_{x})(\varphi^{2^{-n}}_{x})&=(\Pi_{x+h}f(x+h)-\Pi_{x}f(x))(\varphi^{2^{-n}}_{x})\\
	&=\sum_{a \in \mathcal{A}_{\gamma}}\Pi_{x+h} \Big ((f(x+h)-\Gamma_{x+h,x}f(x))_{a}\Big ) (\varphi^{2^{-n}}_{x}) \\
	&=\sum_{a \in \mathcal{A}_{\gamma}}\Pi_{x} \Gamma_{x, x + h} \Big ((f(x+h)-\Gamma_{x+h,x}f(x))_{a}\Big ) (\varphi^{2^{-n}}_{x}) .
	\end{align*}

Using the analytic bounds on $\Pi$ and $\Gamma$, it holds:
	\begin{displaymath}
	|(F_{x+h}-F_{x})(\varphi^{2^{- n}}_{x})|
	\lesssim  \sum_{a \in \mathcal{A}_{\gamma}} 2^{- n a} |f(x+h)-\Gamma_{x+h,x}f(x)|_{a} ,
	\end{displaymath}

\noindent which implies that for $\alpha:=\inf A_{\gamma} \leq a<\gamma$ and $n \in \mathbb{N}$:
	\begin{align*}
	\frac{|(F_{x+h}-F_{x})(\varphi^{2^{- n}}_{x})|}{2^{- n \alpha}(2^{-n} + |h|)^{\gamma-\alpha}}\lesssim \sum_{a \in \mathcal{A}_{\gamma}}2^{- n \left(a-\alpha \right)}\frac{ |f(x+h)-\Gamma_{x+h,x}f(x)|_{a}}{|h|^{\gamma-a}} .
	\end{align*}

Recall that the $\mathcal{D}_{p, q}^{\gamma}$-norm of the modelled distribution $f$ is given by:
	\begin{equation}
		\left\| f \right\|_{\mathcal{D}_{p, q}^{\gamma}} = \sum_{a \in \mathcal{A}_{\gamma}} \left( \left\| \left| f \left(x \right) \right|_{a} \right\|_{L^p \left( x \right)} + \left\| \left\|\frac{ |f(x + h)-\Gamma_{x+h,x}f(x)|_{a}}{|h|^{\gamma-a}} \right\|_{L^{p}(x)} \right\|_{L_h^q \left( h\right)} \right).
	\end{equation}		
	
Thus after integration over $h \in B \left( 0, 1 \right)$:
	\begin{align}
		& \left\| \left\| \left\| \frac{\left( F_{x + h} - F_x \right) \left( \varphi_x^{2^{-n}} \right)}{2^{- n \alpha} \left( 2^{-n} + \left| h \right| \right)^{\gamma - \alpha}} \right\|_{L^{p} \left( x \right)} \right\|_{\ell^{\infty} \left( n \right)} \right\|_{L_h^q \left( h \in B \left( 0, 1 \right) \right)} \\
		& \quad \lesssim \sum_{a \in \mathcal{A}_{\gamma}} \left\| \left\|\frac{ |f(x + h)-\Gamma_{x+h,x}f(x)|_{a}}{|h|^{\gamma-a}} \right\|_{L^{p}(x)} \right\|_{L_h^q \left( h\right)} \lesssim \left\| f \right\|_{\mathcal{D}_{p, q}^{\gamma}} < + \infty ,
	\end{align}

Also, when $h \in B \left( 1, 2 \right)$, using the analytic bound on $\Gamma$ it is straightforward to establish:
	\begin{equation}
		\left\| \frac{\left( F_{x + h} - F_x \right) \left( \varphi_x^{2^{-n}} \right)}{2^{- n \alpha} \left( 2^{-n} + \left| h \right| \right)^{\gamma - \alpha}} \right\|_{L^{p} \left( x \right)} \lesssim \sum_{a \in \mathcal{A}_{\gamma}} \left\| \left| f \left(x \right) \right|_{a} \right\|_{L^p \left( x \right)} \leq \left\| f \right\|_{\mathcal{D}_{p, q}^{\gamma}} ,
	\end{equation}

\noindent whence the property of coherence \eqref{eq:coherence_condition_1}.
Recall that in \cite{MR3684891} it is announced when $q < + \infty$ that for any $\kappa > 0$, $\mathcal{R}^{\rm{HL}} ( f ) \in \mathcal{B}_{p, q}^{\alpha - \kappa}$, while our result yields the seemingly different $\mathcal{R} ( F ) \in \mathcal{B}_{p, \infty}^{\alpha}$.
Our result is actually stronger than those presented in \cite{MR3684891} since $\mathcal{B}_{p, \infty}^{\alpha} \subset \mathcal{B}_{p, 1}^{\alpha - \kappa} \subset \mathcal{B}_{p, q}^{\alpha - \kappa}$.
\end{example}

\begin{example}[Taylor germ]
Let us add one more pedagogical example of germs related to classical Taylor expansions, which have been discussed in \cite[Examples~4.11~\&~5.4]{caravenna2020hairers} in the H\"older case.
Fix $\gamma > 0$ with $\gamma \notin \mathbb{N}$, $p, q \in \left[ 1, + \infty \right]$, and $f \in \mathcal{B}_{p, q}^{\gamma}$.
\autoref{prop:Besov_equivalent_norms} implies that for $0 \leq \left| k \right| < \gamma$, $\partial^k f$ coincides with an $L^{p}$ function (up to a set of Lebesgue measure 0).
Thus we define, for $x, z \in \mathbb{R}^d$ the following germ, called the Taylor germ of $f$:
	\begin{equation} \label{eq:Taylor_germ_of_f_example}
		F_x \left( z \right) \coloneqq \sum\limits_{0 \leq \left| k \right| < \gamma} \partial^k f \left( x \right) \frac{\left( z - x \right)^k}{k !} .
	\end{equation}

Then, $F$ satisfies the properties of coherence and homogeneity in the sense that for all test-functions $\varphi$ with $\int \varphi \neq 0$:
	\begin{equation}
		 \left\| F \right\|_{p, 0, \mathbb{R}^d, \varphi}^{\text{hom}} + \left\| F \right\|_{p, q, 0, \gamma, \mathbb{R}^d, \varphi}^{\text{coh}} < + \infty .
	\end{equation}

Indeed, this corresponds to a particular case of the calculations in \autoref{section:application_young_multiplication} (taking $g \equiv 1$).
Using \autoref{prop:Besov_equivalent_norms}, \autoref{item:prop_equivalent_besov_norm_2} and the uniqueness part in \autoref{corollary:reconstruction_with_coherence_homogeneity}, it is straightforward to observe that $\mathcal{R} \left( F \right) = f$.
\end{example}

\subsection{The problem of Young multiplication in Besov spaces}\label{subsection:young_multiplication}

Now let us discuss the classical problem of multiplying two distributions, provided they belong to suitable Besov spaces, mirroring a similar discussion from \cite[Section~14]{caravenna2020hairers} in the H\"older case.
We will construct such a multiplication as a consequence of our general result \autoref{Prop:Existence_gamma_gq_0}.

The question can be formulated as follows: given $\alpha, \beta, \gamma \in \mathbb{R}$, $p_1$, $p_2$, $p_3$, $q_1$, $q_2$, $q_3 \in \left[ 1, + \infty \right]$, does there exist a continuous bilinear application $\mathcal{M} \colon \mathcal{B}_{p_1, q_1}^{\alpha} \times \mathcal{B}_{p_2, q_2}^{\beta} \to \mathcal{B}_{p_3, q_3}^{\gamma}$ that extends the canonical pointwise multiplication between smooth functions?

Usually, such multiplication maps are constructed with tools from the theory of paraproducts, and it is sometimes claimed in the literature that it is enough to assume:
	\begin{equation}
		\alpha < 0 < \beta, \quad \alpha + \beta > 0, \quad \gamma = \alpha, \quad \frac{1}{p_3} = \frac{1}{p_1} + \frac{1}{p_2} , \quad \frac{1}{q_3} = \frac{1}{q_1} + \frac{1}{q_2} .
	\end{equation}

However, it turns out that only the \enquote{resonant} term of the paraproduct decomposition is well-defined under these conditions.
In fact, the last condition on $q_3$ is incorrect as \cite[Theorem~4.2]{MR1355014} exhibits sequences of smooth functions $f_n, g_n$ such that for any such $\alpha, \beta, p_1, p_2, p_3, q_1, q_2$:
	\begin{equation}
		\begin{dcases}
			\left\| g_n \right\|_{\mathcal{B}_{p_1, q_1}^{\alpha}}, \left\| f_n \right\|_{\mathcal{B}_{p_2, q_2}^{\beta}} & = 1 , \\
			\left\| g_n \cdot f_n \right\|_{\mathcal{B}_{p_3, q}^{\alpha}} & \xrightarrow[n \to + \infty]{} + \infty \text{ when } q < q_1 .
		\end{dcases}
	\end{equation}

In this article, we construct a suitable multiplication map $\mathcal{M}$ under the conditions:
	\begin{equation}
		\alpha < 0 < \beta, \quad \alpha + \beta > 0, \quad \gamma = \alpha, \quad \frac{1}{p_3} = \frac{1}{p_1} + \frac{1}{p_2} , \quad q_3 = q_1 .
		\label{eq:conditions_parameters_Young_multiplication}
	\end{equation}

That is, we build a multiplication:
	\begin{equation}
		\mathcal{M} \colon \mathcal{B}_{p_1, q_1}^{\alpha} \times \mathcal{B}_{p_2, q_2}^{\beta} \to \mathcal{B}_{\left( p_1^{-1} + p_2^{-1} \right)^{-1}, q_1}^{\alpha} .
	\end{equation}

Note that it is known that such a map $\mathcal{M}$ can be built with paraproducts, see \cite[Theorem~6.6]{MR1355014} or \cite[Corollary~2.1.35]{Martin2018Refinements}.
However, our construction does not require paraproducts and relies instead on \autoref{Prop:Existence_gamma_gq_0} below.
We will provide a proof in \autoref{section:application_young_multiplication}, but let us outline the strategy here.

Recalling the embedding $\mathcal{B}_{p_2, q_2}^{\beta} \subset \mathcal{B}_{p_2, q_2}^{\beta - \epsilon}$ for any $\epsilon > 0$, we can assume without loss of generality that $\beta \notin \mathbb{N}$.

Thus, we shall fix $\alpha < 0, \beta > 0$ with $\alpha + \beta > 0$ and $\beta \notin \mathbb{N}$, as well as $p_1, p_2, q_1, q_2, p, q \in \left[ 1, + \infty \right]$ with $\frac{1}{p} = \frac{1}{p_1} + \frac{1}{p_2}$, $\frac{1}{q} = \frac{1}{q_1} + \frac{1}{q_2}$.
Fix distributions $g \in \mathcal{B}_{p_1, q_1}^{\alpha}$, $f \in \mathcal{B}_{p_2, q_2}^{\beta}$.
Note that our conventions for the sign of $\alpha$ and $\beta$ are interverted with those of \cite{caravenna2020hairers}.

Since $f \in \mathcal{B}_{p_2, q_2}^{\beta}$ with $\beta > 0$ and $\beta \notin \mathbb{N}$, we know from \autoref{prop:Besov_equivalent_norms} that for $0 \leq \left| k \right| < \beta$, $\partial^k f$ coincides with an $L^{p_2}$ function (up to a set of Lebesgue measure 0), so that we can define, for $x, z \in \mathbb{R}^d$:
	\begin{equation} \label{eq:Taylor_germ_of_f}
		F_x \left( z \right) \coloneqq \sum\limits_{0 \leq \left| k \right| < \beta} \partial^k f \left( x \right) \frac{\left( z - x \right)^k}{k !} .
	\end{equation}

This is the Taylor germ of $f$.
Now we define a germ $P$ by setting for $x \in \mathbb{R}^d$ and $\varphi \in \mathcal{D} ( \mathbb{R}^d )$:
	\begin{equation} \label{eq:germ_for_multiplication}
		P_x \left( \varphi \right) \coloneqq g \left( \varphi F_x \right) .
	\end{equation}

Recall that this is the same germ as considered in \cite{caravenna2020hairers} in the case of Young multiplication for H\"older distributions.
Note also that here, $P_x$ is only correctly defined for $x$ away from a null set, which is not a problem since all the objects required for the reconstruction are defined by integration over $x$.

In \autoref{section:application_young_multiplication} we will prove that the germ $P$ satisfies the hypotheses of the more general Reconstruction \autoref{Prop:Existence_gamma_gq_0} below, and that $P$ admits a unique reconstruction $\mathcal{R} \left( P \right)$.
Therefore the following result will follow by setting $\mathcal{M} \left( g, f \right) \coloneqq \mathcal{R} \left( P \right)$:
\begin{theorem}[Young multiplication in Besov spaces]\label{prop:Young_multiplication_besov}
Let $p_1, p_2, q_1, q_2 \in \left[ 1, + \infty \right]$ and let $p, q \in \left[ 1, + \infty \right]$ be defined by $\frac{1}{p} = \frac{1}{p_1} + \frac{1}{p_2}$, $\frac{1}{q} = \frac{1}{q_1} + \frac{1}{q_2}$.
Let $\alpha, \beta \in \mathbb{R}$ be such that $\alpha < 0 < \beta$, $\alpha + \beta > 0$.

Then there exists a bilinear continuous map $\mathcal{M} \colon \mathcal{B}_{p_1, q_1}^{\alpha} \times \mathcal{B}_{p_2, q_2}^{\beta} \to \mathcal{B}_{p, q_1}^{\alpha}$ that extends the usual product, i.e.\ when $g \in \mathcal{B}_{p_1, q_1}^{\alpha}$ and $f \in \mathcal{C}^{\infty}$, $\mathcal{M} \left( g, f \right) = g \cdot f$, where the product in the right-hand side is understood as the product of a distribution against a smooth function.

Furthermore, when $\beta \notin \mathbb{N}$, our map $\mathcal{M}$ is characterised by the following property: for any $r \in \mathbb{N}$ with $r > - \alpha$:
	\begin{equation} \label{eq:reconstruction_bound_Young_multiplication}
		\left\| \left\| \sup\limits_{\psi \in \mathcal{B}^r} \left| \frac{\left( \mathcal{M} \left( g, f \right) - g \cdot F_x \right) \left( \psi_x^{\lambda} \right)}{\lambda^{\alpha + \beta}} \right| \right\|_{L^p \left( x \right)} \right\|_{L_{\lambda}^q \left( \lambda \right)} < + \infty .
	\end{equation}
\end{theorem}

\begin{remark}
A similar bound as \eqref{eq:reconstruction_bound_Young_multiplication} is established in \cite[Equation~(3.1)]{liu2021sobolev} in the case $p_1 = p_2 = q_1 = q_2$.
\end{remark}

\section{A general Besov reconstruction theorem}
\label{section:proof_of_main_result}

Now let us turn to the statement and proof of our most general reconstruction result, \autoref{Prop:Existence_gamma_gq_0} below. 

\subsection{Statement of the result}

Let us introduce the following notations.

\begin{definition}[Besov reconstruction of a Germ]\label{def:Besov_reconstruction}
Let $F$ be a germ, $p, q \in \left[ 1, + \infty \right]$, and $r \in \mathbb{N}$.
Let $k \colon \left( 0, 1 \right] \to \mathbb{R}_+$ be a function (which we will call a scaling function).
Let $K \subset \mathbb{R}^d$.
We say that a distribution $\mathcal{R} \left( F \right) \in \mathcal{D}^{\prime} ( \bar{K}_1 )$ is a \emph{$k, p, q$-reconstruction of $F$ on $K$} if the following estimate, called the \emph{reconstruction bound}, holds for any test-function $\psi \in \mathcal{D} ( B ( 0, 1 ) )$:
	\begin{equation}
		\left\| \left\| \frac{\left( \mathcal{R} \left( F \right) - F_x \right) \left( \psi_x^{\lambda} \right)}{k \left( \lambda \right)} \right\|_{L^{p} \left( x \in K \right)} \right\|_{L_{\lambda}^q(\lambda)} < + \infty .
		\label{eq:definition_reconstruction_bound_in_lambda}
	\end{equation}
	
	We say that the reconstruction $\mathcal{R} \left( F \right)$ is \emph{$r$-uniform} if the reconstruction bound is uniform in $\psi$ in the following sense:
	\begin{equation}
		\left\| \left\| \sup\limits_{\psi \in \mathscr{B}^{r}} \left| \frac{\left( \mathcal{R} \left( F \right) - F_x \right) \left( \psi_x^{\lambda} \right)}{k \left( \lambda \right)} \right| \right\|_{L^{p} \left( x \in K \right)} \right\|_{L^q_\lambda(\lambda)} < + \infty .
		\label{eq:definition_reconstruction_bound_in_lambda_uniform_in_psi}
	\end{equation}

	Finally, if $\gamma \in \mathbb{R}$, we say that a distribution is a \emph{$\gamma, p, q$-reconstruction of $F$ on $K$} if it is a $k_{\gamma, q, \epsilon}, p, q$-reconstruction of $F$ on $K$, where $k$ is defined by \eqref{eq:definition_scaling_function}.
\end{definition}

The reconstruction bound quantifies how close $\mathcal{R} \left( F \right)$ is to $F_x$ locally in space (the $x$ variable) and scale (the $\lambda$ variable), in a way similar to the definition of Besov spaces, see \autoref{appendix_besov_spaces}.

It is interesting to note that this definition already guarantees some properties of the reconstruction:
\begin{proposition}\label{Prop:Uniqueness}
In the context of the previous definition, let $k_1, k_2 \colon \left( 0, 1 \right] \to \mathbb{R}_+$ be two scaling functions.
	\begin{enumerate}[ref =\emph{(\arabic*)}]
		\item\label{item:scaling_prop_1} Assume there exists $C > 0$ such that $k_1 \leq C k_2$ pointwise, then a $k_1, p, q$-reconstruction of $F$ is also a $k_2, p, q$-reconstruction of $F$.
		\item\label{item:scaling_prop_2} Assume there exists $C > 0$ such that for all $\lambda \in \left( 0, 1 \right]$, $C^{-1} k \left( \lambda \right) \leq k \left( 2^{ \left\lfloor \log_2 \left( \lambda \right) \right\rfloor} \right) \leq C k \left( \lambda \right)$. Then a distribution $\mathcal{R} \left( F \right) \in \mathcal{D}^{\prime} ( \bar{K}_1 )$ is a $r$-uniform $k, p, q$-reconstruction of $F$ on $K$ if and only if:
			\begin{equation}
			\left\| \left\| \sup\limits_{\psi \in \mathscr{B}^{r}} \left| \frac{\left( \mathcal{R} \left( F \right) - F_x \right) \left( \psi_x^{2^{-n}} \right)}{k \left( 2^{-n} \right)} \right| \right\|_{L^{p} \left( x \in K \right)} \right\|_{\ell^q \left( n \right)} < + \infty .
			\label{eq:definition_reconstruction_bound_in_n}
		\end{equation}
		\item\label{item:scaling_prop_3} Assume that $k_1, k_2$ are two scaling function such that \autoref{item:scaling_prop_2} just above applies, and $k_1, k_2 \left( \lambda \right) \to_{\lambda \to 0} 0$. Then any $k_1, p, q$-reconstruction of a germ $F$ on all $K \subset \mathbb{R}^d$ must coincide with any $k_2, p, q$-reconstruction of $F$ on all $K \subset \mathbb{R}^d$. In particular, when $k \left( \lambda \right) \to_{\lambda \to 0} 0$, a germ $F$ can have at most one $k, p, q$-reconstruction.
	\end{enumerate}
\end{proposition}

\begin{proof}
\ref{item:scaling_prop_1} and \ref{item:scaling_prop_2} are elementary to check.
Now let us tackle \ref{item:scaling_prop_3}.
Denote $T$ the difference of the two reconstructions.
The embedding $\ell^q\subset \ell^\infty$, H\"older's inequality, and the triangle inequality, yield $\lim_{n \to + \infty} \int_{K} | T ( \varphi_x^{\epsilon_n} ) | dx = 0$ for any compact $K \subset \mathbb{R}^d$ and any test-function $\varphi$.
Now we fix a test-function $\varphi$ such that $\int \varphi = 1$.
Let $\eta \in \mathcal{D} ( \mathbb{R}^d )$ and let us show that $T ( \eta ) = 0$.
By mollification, $T ( \eta ) = \lim_{n \to + \infty} T_n ( \eta )$ where $T_n ( \eta ) \coloneqq \int_{\supp \left( \eta \right)} T ( \varphi_x^{\epsilon_n} ) \eta ( x ) d x$. Now $| T_n ( \eta ) | \leq \left\| \eta \right\|_{\infty} \int_{\supp \left( \eta \right)} | T ( \varphi_x^{\epsilon_n} ) | dx = o_{n \to + \infty} ( 1 )$, hence $T ( \eta ) = 0$ as announced, so that the two reconstructions coincide.
\end{proof}

\begin{notation}
If $F$ is a germ, $p \in \left[ 1, + \infty \right]$, $\alpha, \beta, \gamma \in \mathbb{R}$, $K \subset \mathbb{R}^d$, and $\varphi \in \mathcal{D} ( \mathbb{R}^d )$ is a test-function, we denote:
	\begin{equation}
		\begin{dcases}
			f_K \left( n, h \right) \coloneqq f_{F, \varphi, \alpha, \gamma, p, K} \left( n, h \right) & \coloneqq \left\| \frac{\left( F_{x + h} - F_x \right) \left( \varphi_x^{2^{-n}} \right)}{2^{- n \alpha} \left( 2^{-n} + \left| h \right| \right)^{\gamma - \alpha}} \right\|_{L^{p} \left( x \in \bar{K}_2 \right)} , \\
			g_K \left( n \right) \coloneqq g_{F, \varphi, \beta, p, K} \left( n \right) & \coloneqq \left\| \frac{F_x \left( \varphi_x^{2^{-n}} \right)}{2^{- n \beta}} \right\|_{L^{p} \left( x \in \bar{K}_2 \right)} .
		\end{dcases}	
		\label{eq:definition_of_f_and_g_in_the_proposition}
	\end{equation}

We will usually drop the subscripts when the dependence in $F, \varphi, \alpha, \beta, \gamma, p$ is clear from the context.
\end{notation}

The following sequences, corresponding to \emph{averaged} versions of $f$, will play an important role in our calculations.
Note that these quantities will appear naturally in our calculations, but we do not really have an interpetation of what they represent.
\begin{notation}
For any function $f \colon \mathbb{N} \times \mathbb{R}^d \to \mathbb{R}_+$ and real $c \in \mathbb{R}$, set for $n \in \mathbb{N}$:
	\begin{equation}
		\begin{dcases}
			m^{\left( 1 \right)}_{f} \left( n \right) & \coloneqq \int_{h \in B \left( 0, 2^{-n + 1} \right)} 2^{n d} f \left( n, h \right) d h , \\
			m^{\left( 2 \right)}_{c, f} \left( n \right) & \coloneqq \sum\limits_{k = n}^{+ \infty} \int_{h \in B \left( 0, 2^{-k} \right)} 2^{- \left( k - n \right) c + k d} f \left( k, h \right) d h , \\
			m^{\left( 3 \right)}_{c, f} \left( n \right) & \coloneqq \sum\limits_{k = n}^{+ \infty} \int_{h \in B \left( 0, 2^{-n + 1} \right)} 2^{- \left( k - n \right) c + n d} f \left( k, h \right) d h , \\
			m^{\left( 4 \right)}_{c, f} \left( n \right) & \coloneqq \sum\limits_{k = 0}^{n - 1} \int_{h \in B \left( 0, 2^{-k + 1} \right)} 2^{- \left( k - n \right) c + k d} f \left( k, h \right) d h .
		\end{dcases}
		\label{eq:definition_of_m1_m2_m3}
	\end{equation}

Once again, we will drop the subscripts when the dependence in the parameters is clear from context.
\end{notation}

We will establish the following general version of the Besov reconstruction theorem.
\begin{theorem}[Besov reconstruction, general case] \label{Prop:Existence_gamma_gq_0}
Let $F$ be a germ, $p, q \in \left[ 1, + \infty \right]$, and $\alpha, \beta, \gamma \in \mathbb{R}$ be such that $\alpha \leq \gamma$.
Let $r \in \mathbb{N}$ be an integer such that $r > - \beta$.
Let $K \subset \mathbb{R}^d$ and $\hat{\varphi} \in \mathcal{D} ( B ( 0, 1 / 2 ) )$ be a test-function such that $\int_{} \hat{\varphi} = 1$ and $\int x^k \hat{\varphi} \left( x \right) d x = 0$ for $1 \leq \left| k \right| \leq r - 1$.
For simplicity of notation, denote:
\begin{alignat}{2}\label{eq:convention_for_definitions}
	\begin{cases} f_K &= f_{F, \hat{\varphi}, \alpha, \gamma, p, K} , \\
			g_K &= g_{F, \hat{\varphi}, \beta, p, K} , \\
			m_K^{\left( 1 \right)} &= m^{\left( 1 \right)}_{f_K} , 							\end{cases} 
	&& \quad \quad \quad \quad \begin{cases} m_K^{\left( 2 \right)} &= m^{\left( 2 \right)}_{\gamma, f_K} , \\
			m_K^{\left( 3 \right)} &= m^{\left( 3 \right)}_{\alpha + r, f_K} , \\
			m_K^{\left( 4 \right)} &= m^{\left( 4 \right)}_{\gamma, f_K} . 
	\end{cases}
\end{alignat}

For any $q_1 \in \left[ 1, + \infty \right]$, we denote:
	\begin{equation}
		\left\| F \right\|_{\mathcal{G}_{p, q, q_1, K, \varphi}^{\alpha, \beta, \gamma}} \coloneqq
		\begin{dcases}
			\left\| g_K \right\|_{\ell^{q_1}} + \left\| m_K^{\left( 1 \right)} \right\|_{\ell^q} + \left\| m_K^{\left( 2 \right)} \right\|_{\ell^q} + \left\| m_K^{\left( 3 \right)} \right\|_{\ell^q} & \text{ if } \gamma > 0 , \\
				\left\| g_K \right\|_{\ell^{q_1}} + \left\| m_K^{\left( 3 \right)} \right\|_{\ell^q} + \left\| m_K^{\left( 4 \right)} \right\|_{\ell^q} & \text{ if } \gamma < 0 , \\
				\left\| g_K \right\|_{\ell^{q_1}} + \left\| \tilde{m}_K^{\left( 3 \right)} \right\|_{\ell^q} + \left\| \tilde{m}_K^{\left( 4 \right)} \right\|_{\ell^q} & \text{ if } \gamma = 0 ,
		\end{dcases}
	\end{equation}

\noindent where in the case $\gamma = 0$, we define, for $n \in \mathbb{N}$ and $i \in \left\lbrace 3, 4 \right\rbrace$, $\tilde{m}_K^{\left( i \right)} \left( n \right) \coloneqq m^{\left( i \right)} \left( n \right)/k \left( 2^{-n} \right)$, where $k$ is any scaling function such that $\left( 1/{k \left( 2^{- n} \right)} \right)_{n \in \mathbb{N}} \in \ell^q$.

Assume that:
	\begin{equation}\label{eq:assumption_of_main_theorem}
		\left\| F \right\|_{\mathcal{G}_{p, q, q_1, K, \varphi}^{\alpha, \beta, \gamma}} < + \infty .
	\end{equation}

Then there exists a $k, p, q$-reconstruction of $F$ on $K$, noted $\mathcal{R} \left( F \right)$ or $\mathcal{R}_K \left( F \right)$, that is also $r$-uniform.

Furthermore:
	\begin{enumerate}
		\item (Global version) if \eqref{eq:assumption_of_main_theorem} holds for $K = \mathbb{R}^d$, then:
		\begin{enumerate}
			\item Such an $\mathcal{R} \left( F \right)$ is unique when $\gamma > 0$ but not when $\gamma \leq 0$.
			\item If $\beta \wedge \gamma > 0$, then $\mathcal{R} ( F ) = 0$.
			\item If $\beta \wedge \gamma = 0$, then for all $\kappa > 0$, $\mathcal{R} \left( F \right) \in \mathcal{B}_{p, 1}^{- \kappa}$.
			\item If $\beta \wedge \gamma < 0$, then $\mathcal{R} \left( F \right) \in \mathcal{B}_{p, q_1 \vee q}^{\beta \wedge \gamma}$.
		\end{enumerate}
	
	\item (Local version) if \eqref{eq:assumption_of_main_theorem} holds for all $K \subset \mathbb{R}^d$, then there exists a global distribution $\mathcal{R} \left( F \right) \in \mathcal{D}^{\prime} ( \mathbb{R}^d )$ that is a $r$-uniform $k, p, q$-reconstruction of $F$ on all $K \subset \mathbb{R}^d$ and:
		\begin{enumerate}
			\item Such an $\mathcal{R} \left( F \right)$ is unique when $\gamma > 0$ but not when $\gamma \leq 0$.
			\item If $\beta \wedge \gamma > 0$, then $\mathcal{R} ( F ) = 0$.
			\item If $\beta \wedge \gamma = 0$, then for all $\kappa > 0$, $\mathcal{R} \left( F \right) \in \mathcal{B}_{p, 1, \mathrm{loc}}^{- \kappa}$.
			\item If $\beta \wedge \gamma < 0$, then $\mathcal{R} \left( F \right) \in \mathcal{B}_{p, q_1 \vee q, \mathrm{loc}}^{\beta \wedge \gamma}$.
		\end{enumerate}
	\end{enumerate}

The reconstruction map is continuous in the following sense: let $\mathcal{B}$ denote $\mathcal{B}_{p, 1, K}^{- \kappa}$ if $\beta \wedge \gamma = 0$ and $\mathcal{B}_{p, q_1 \vee q, K}^{\beta \wedge \gamma}$ if $\beta \wedge \gamma < 0$, then:
	\begin{equation}
		\left\| \mathcal{R} \left( F \right) \right\|_{\mathcal{B}} \lesssim \left\| F \right\|_{\mathcal{G}_{p, q, q_1, \bar{K}_2, \varphi}^{\alpha, \beta, \gamma}} .
	\end{equation}
\end{theorem}

\begin{remark}
The cases $\gamma \neq 0$ could also be slightly modified to consider general scaling functions $k$ as in the case of $\gamma = 0$.
\end{remark}

Let us turn to the proof of \autoref{Prop:Existence_gamma_gq_0}.
In the remainder of this section, we consider a germ $F$, reals $p, q \in \left[ 1, + \infty \right]$, $\alpha, \beta, \gamma \in \mathbb{R}$ such that $\alpha \leq \gamma$, an integer $r \in \mathbb{N}$, $K \subset \mathbb{R}^d$, and a single test-function $\hat{\varphi} \in \mathcal{D} ( B ( 0, 1 / 2 ) )$ such that $\int_{} \hat{\varphi} = 1$ and $\int x^k \hat{\varphi} \left( x \right) d x = 0$ for $1 \leq \left| k \right| \leq r - 1$.
We denote $f$, $g$, $m^{\left( 1 \right)}$, $m^{\left( 2 \right)}$, $m^{\left( 3 \right)}$, $m^{\left( 4 \right)}$ to be the functions defined in \eqref{eq:convention_for_definitions}.

We break up the proof into several sections.

\subsection{Uniqueness of reconstruction}

Recall from \autoref{Prop:Uniqueness}, \autoref{item:scaling_prop_3}, that when $\gamma > 0$, the reconstruction, if it exists, is unique. 
Nevertheless, when $\gamma \leq 0$, the reconstruction is not unique in general.
We now focus on the existence of $\mathcal{R} \left( F \right)$.

\subsection{Existence for \texorpdfstring{$\gamma>0$}{gamma > 0}}

We now construct a reconstruction in the case $\gamma > 0$.
First, let us recall the strategy of \cite{caravenna2020hairers}. 

In order to establish the existence of the reconstruction, we proceed by mollification. Recall that if $\rho \in \mathcal{D} ( \mathbb{R}^d )$ is any test-function such that $\int \rho = 1$, and if $\xi \in \mathcal{D}^{\prime} ( \mathbb{R}^d )$ is any distribution, then we have an approximation of $\xi$ provided for $\psi \in \mathcal{D} ( \mathbb{R}^d )$ by: 
\begin{displaymath}
	\xi(\psi)=\lim_{n\rightarrow \infty}\xi \left( \rho^{2^{-n}}*\psi \right) =\lim_{n\rightarrow \infty}\int_{\R^d}\xi \left( \rho_{z}^{2^{-n}} \right) \psi(z)\,dz .
\end{displaymath} 

This yields a natural candidate for the reconstruction of the germ $F$, as we would like to set by analogy:
	\begin{equation}
		\mathcal{R} \left( F \right) \left( \psi \right) \stackrel{\text{?}}{\coloneqq} \lim_{n\rightarrow \infty}\int_{\R^d}F_z \left( \rho_{z}^{2^{-n}} \right) \psi(z)\,dz .
	\end{equation}

Of course, the convergence of such a sequence is a priori far from obvious.
However, as our goal is an existence result for $\mathcal{R} ( F )$, it is enough for us to exhibit just one choice of $\rho$ such that this sequence converges. 
For this purpose, we follow the strategy of \cite{caravenna2020hairers} i.e.\ construct a specific mollifier $\rho$ from the single test-function $\hat{\varphi}$ provided by the assumption of the theorem. 

Explicitly, define as in \cite{caravenna2020hairers}:
\begin{equation}\label{EQN:MOLLIFIER}
	\rho:=\hat{\varphi}^{2}*\hat{\varphi}.
\end{equation}

The motivation behind this choice is that it allows us to rewrite the difference $\rho^{1/2}-\rho$ as a convolution, paving the way for a nice dyadic decomposition. Explicitly, set $\check{\varphi}:=\hat{\varphi}^{1/2}-\hat{\varphi}^{2}$, then:
\begin{displaymath}
	\rho^{1/2}-\rho=\hat{\varphi}*\check{\varphi} .
\end{displaymath}

This directly implies that for $n \in \mathbb{N}$,
\begin{equation}\label{EQN:DIFF_CONV}
	\rho^{2^{-(n+1)}}-\rho^{2^{-n}}=\hat{\varphi}^{2^{-n}}*\check{\varphi}^{2^{-n}} .
\end{equation}

Notice that by assumption on $\hat{\varphi}$, it holds that $\supp \left( \check{\varphi} \right) \subset B \left( 0, 1 \right)$, and that $\check{\varphi}$ cancels all polynomials of degree less that $r - 1$: $\int x^k \hat{\varphi} \left( x \right) d x = 0$ for $0 \leq \left| k \right| \leq r - 1$.

Now we can exploit the dyadic structure of our mollifier: for all $n_0, n \in \mathbb{N}$, 
	\begin{equation}
		\rho^{2^{-n}} = \rho^{2^{- n_0}} + \sum\limits_{k = n_0}^{n - 1} \hat{\varphi}^{2^{-k}}*\check{\varphi}^{2^{-k}} .
	\end{equation}

This gives us a natural definition of approximating our reconstruction. 

\begin{definition}\label{DEF:APPROX}
	Let $F=(F_{x})_{x\in \R^d}$ be a germ.
	For simplicity of notation, denote $\epsilon_n \coloneqq 2^{-n}$.
	We define a sequence of approximating distributions $\mathcal{R}_{n} \left( F \right) \in \mathcal{D}'(\R^d)$ by setting, for $n \in \mathbb{N}$, $\psi \in \mathcal{D} ( \mathbb{R}^d )$, and any $n_0 \in \mathbb{N}$:
	\begin{align}
		\mathcal{R}_{n} \left( F \right)(\psi) & \coloneqq \int_{\R^d}F_{z}\left (\rho_{z}^{\epsilon_n}\right )\psi(z)\,dz \nonumber \\
		& = \int_{\R^d} F_{z}\left (\hat{\varphi}^{\epsilon_{n_0}}*\hat{\varphi}_z^{2\epsilon_{n_0}}\right )\psi(z)\,dz + \sum\limits_{k = n_0}^{n - 1} \int_{\R^d}F_{z}\left (\hat{\varphi}^{\epsilon_k}*\check{\varphi}_z^{\epsilon_k}\right )\psi(z)\,dz .\nonumber
	\end{align}
	where $\rho$ is defined as in \eqref{EQN:MOLLIFIER}.
	Note that as explained above, this definition does not depend on the choice of $n_0 \in \mathbb{N}$.
	If the sequence converges, then we denote for $\psi \in \mathcal{D} ( \mathbb{R}^d )$:
		\begin{equation}
			\mathcal{R}(F)(\psi):=\lim_{n\rightarrow\infty}\mathcal{R}_n \left( F \right) (\psi) .
		\end{equation}
\end{definition}

Now, we want to establish whether this limit $\lim_{n\rightarrow \infty} \mathcal{R}_{n} \left( F \right)$ exists.
For this, we shall pursue even further the decomposition of $\mathcal{R}_{n} \left( F \right)$.
Recall that for any distribution $\xi \in \mathcal{D}^{\prime} ( \mathbb{R}^d )$ and any two test-functions $\eta, \tilde{\eta} \in \mathcal{D} ( \mathbb{R}^d )$, 
	\begin{equation}
		\xi \left( \eta * \tilde{\eta} \right) = \int_{\mathbb{R}^d} \xi \left( \eta_x \right) \tilde{\eta} \left( x \right) d x .
		\label{eq:distribution_against_convolution}
	\end{equation}

From \autoref{DEF:APPROX}, it follows that the existence of $\mathcal{R} ( F )$ is implied by the absolute convergence of the series $\sum_{k} u_k$, where we set for $k \in \mathbb{N}$:
	\begin{align}
		u_k & \coloneqq \int_{\R^d}F_{z}\left (\hat{\varphi}^{\epsilon_k}*\check{\varphi}_z^{\epsilon_k}\right )\psi(z)\,dz \\
		& = \int_{\R^d}\int_{\R^d}F_{z}\left(\hat{\varphi}^{\epsilon_k}_{x}\right)\check{\varphi}^{\epsilon_k}(x-z)\psi(z)\,dx\,dz .
	\end{align}

Writing $F_{z}=F_{x}+(F_{z}-F_{x})$, we decompose $u_k = u_k^{\prime} + u_k^{\prime \prime}$, where:
\begin{align}
	u'_k& \coloneqq \int_{\R^d}\int_{\R^d}F_{x}\left(\hat{\varphi}^{\epsilon_k}_{x}\right)\check{\varphi}^{\epsilon_k}(x-z)\psi(z)\,dx\,dz\\ 
	&=\int_{\R^d}F_{x}\left(\hat{\varphi}^{\epsilon_k}_{x}\right)\left (\check{\varphi}^{\epsilon_k}*\psi\right )(x)\,dx ,
\end{align}

\noindent and:
\begin{align}	\label{EQN:uii}
	u''_k& \coloneqq \int_{\R^d}\int_{\R^d}(F_{z}-F_{x})\left(\hat{\varphi}^{\epsilon_k}_{x}\right)\check{\varphi}^{\epsilon_k}(x-z)\psi(z)\,dx\,dz\\
	&=\int_{\R^d}\int_{\R^d}(F_{x+h}-F_{x})\left(\hat{\varphi}^{\epsilon_k}_{x}\right)\check{\varphi}^{\epsilon_k}(-h)\psi(x+h)\,dx\,dh . \nonumber
\end{align}

Hence, the existence of the reconstruction $\mathcal{R} ( F )$ can be determined by the absolute convergence of
\begin{displaymath}
	\sum_{k}u'_k, \quad \sum_{k}u''_k. 
\end{displaymath}

The following lemma from \cite{caravenna2020hairers} will also be useful for us here. For completeness, we concisely recall its proof. Note that \cite[Lemma 9.2]{caravenna2020hairers} states the bound in the $L^1$ case but its proof actually treats the $L^{\infty}$ case. 

\begin{lemma}\cite[Lemma 9.2]{caravenna2020hairers}\label{LEMMA:CONV_INEQ}
	Assume that $\check{\varphi} \in \mathcal{D} ( \mathbb{R}^d )$ is a test-function that cancels polynomials of degree $r- 1 \in \mathbb{N}$ i.e.\ for all $0 \leq k \leq r - 1$, $\int_{\mathbb{R}^d} x^k \check{\varphi} \left( x \right) d x = 0$. 
	Then for any test-function $\eta \in \mathcal{D}( \mathbb{R}^d )$ and $\lambda >0$:
	\begin{displaymath}
		\|\check{\varphi}^{\lambda}*\eta\|_{L^\infty}\leq \|\eta\|_{C^r}\|\check{\varphi}\|_{L^1}\lambda^{r}.
	\end{displaymath}
\end{lemma}

\begin{proof}
	We fix $y\in \mathbb{R}^d$ and denote $p_{y}(\cdot) \coloneqq \sum_{|k|\leq r-1}\frac{\partial^k\eta(y)}{k!}(\cdot - y)^k$ to be the Taylor polynomial of $\eta$ of order $r-1$ based at $y$, so that $|\eta(z)-p_{y}(z)|\leq \|\eta\|_{C^r}|z-y|^r$. As noted above, $\check{\varphi}$ cancels polynomials of degree $r - 1$, so that $\int_{\R^d}\check{\varphi}^\lambda(y-z)p_{y}(z) dz=0$, thus:
	\begin{displaymath}
		(\check{\varphi}^{\lambda}*\eta)(y)=\int_{\R^d}\check{\varphi}^\lambda(y-z)\{\eta(z)-p_{y}(z)\}\,dz.
	\end{displaymath}
	Hence,
	\begin{align*}
		|(\check{\varphi}^{\lambda}*\eta)(y)|
		&\leq \int_{\R^d}|\check{\varphi}^\lambda(y-z)||\eta(z)-p_{y}(z)|\,dz\\
		&\leq \|\eta\|_{C^r}\int_{\R^d}|\check{\varphi}^\lambda(y-z)||z-y|^r\,dz\\
		&=	\lambda^r\|\eta\|_{C^r}\|\check{\varphi}\|_{L^1},
	\end{align*}
	which gives us our result. 
\end{proof}

\begin{proposition}[Convergence of approximating distributions]
\label{prop:convergence_approximating_distributions}
In the setting of this section:
	\begin{enumerate}[ref =\textit{(\arabic*)}]
		\item\label{item:convergence_1} Suppose that $\beta + r > 0$ and that $g_K \in \ell^{\infty}$. Then $\sum_{k}u'_k$ is absolutely convergent as soon as $\supp \left( \psi \right) \subset \bar{K}_1$. 
		\item\label{item:convergence_2} Suppose that $\gamma > 0$, $\gamma \geq \alpha$, and $m_K^{(1)} \in \ell^q$. Then $\sum_{k}u_k''$ is absolutely convergent as soon as $\supp \left( \psi \right) \subset \bar{K}_1$.
		\item\label{item:convergence_3} Suppose that assumptions \ref{item:convergence_1} and \ref{item:convergence_2} just above apply. Then $\mathcal{R} ( F ) \in \mathcal{D}^{\prime} ( \bar{K}_1 )$ is a distribution of order $r$.
	\end{enumerate}

\end{proposition}
\begin{proof}
We start with \ref{item:convergence_1}.
	Denote $\tilde{p}$ to be the H\"older conjugate of $p$. By H\"older's inequality, and since $\supp \left( \check{\varphi}^{\epsilon_k}*\psi \right) \subset \bar{K}_2$
	\begin{align*}
		|u'_k|&\leq \int_{\bar{K}_2}\left|F_{x}\left(\hat{\varphi}^{\epsilon_k}_{x}\right)\left (\check{\varphi}^{\epsilon_k}*\psi\right )(x)\right |\,dx\\
		&\leq \left \|F_{x}\left(\hat{\varphi}^{\epsilon_k}_{x}\right)\right\|_{L^p(x \in \bar{K}_2)}\left\| \check{\varphi}^{\epsilon_k}*\psi\right \|_{L^{\tilde{p}}} .
	\end{align*}

	Applying \autoref{LEMMA:CONV_INEQ}, we obtain: 
	\begin{align}
		\left\| \check{\varphi}^{\epsilon_k}*\psi\right \|_{L^{\tilde{p}}} \lesssim \|\psi\|_{C^{r}}\|\check{\varphi}\|_{L^1}2^{-kr} . 
	\end{align}
	Thus, recalling the definition of $g_K$, and since $\beta + r > 0$,
	\begin{align*}
		\sum_{k}|u'_k|&\lesssim \|\psi\|_{C^r} \|\check{\varphi}\|_{L^1} \sum_{k} \left( \frac{\left\|F_{x}\left(\hat{\varphi}^{\epsilon_k}_{x}\right)\right\|_{L^p(x \in \bar{K}_2)}}{2^{-k\beta}} \right) 2^{-k(\beta+r)}\\
		&\lesssim \|\psi\|_{C^r} \|\check{\varphi}\|_{L^1} \|g_K\|_{\ell^\infty}\sum_{k}2^{-k(\beta+r)}<\infty. 
	\end{align*}
This yields the announced result.

Now we prove \ref{item:convergence_2}.
By definition:
	\begin{align}
		|u''_k|&\leq \int_{\R^d}\int_{\R^d}\left |(F_{x+h}-F_{x})\left(\hat{\varphi}^{\epsilon_k}_{x}\right)\check{\varphi}^{\epsilon_k}(-h)\psi(x+h)\right |\,dx\,dh
	\end{align}

	Because of the supports of $\check{\varphi}$	and $\psi$, we have that $h$ runs over $B \left( 0, \epsilon_k \right)$ and $x$ runs over $\bar{K}_2$.
	If we denote $\tilde{p}$ to be the H\"older conjugate of $p$, then by H\"older's inequality, with respect to the $x$ variable, 
	\begin{align}
		|u''_k|&\leq \int_{B \left( 0, \epsilon_k \right)}\left \|(F_{x+h}-F_{x})\left(\hat{\varphi}^{\epsilon_k}_{x}\right)\right \|_{L^p(x \in \bar{K}_2)}\left \|\psi(x+h)\right \|_{L^{\tilde{p}}(x)}|\check{\varphi}^{\epsilon_k}(-h)|\,dh. \nonumber
	\end{align}	
	
	By substitution, $\left \|\psi(x+h)\right \|_{L^{\tilde{p}}(x)}=\|\psi\|_{L^{\tilde{p}}}$.
	Also, $|\check{\varphi}^{\epsilon_k}(-h)|\leq 2^{d k}\|\check{\varphi}\|_{L^\infty}$, so that:
	\begin{align}	
		|u_k''|&\leq 2^{dk}\|\check{\varphi}\|_{L^\infty} \|\psi\|_{L^{\tilde{p}}}\int_{B \left( 0, \epsilon_k \right)}\left \|(F_{x+h}-F_{x})\left(\hat{\varphi}^{\epsilon_k}_{x}\right)\right \|_{L^p(x)}\,dh. 
	\end{align}
	Recall that by definition of $f_K \left( k, h \right)$:
		\begin{align}
			\left\|(F_{x+h}-F_{x})\left(\hat{\varphi}^{\epsilon_k}_{x}\right)\right\|_{L^p(x\in \bar{K}_2)} &= f_K \left( k, h \right) 2^{- k \alpha} (2^{-k}+|h|)^{\gamma-\alpha} ,
		\end{align}
	\noindent and since here $\left| h \right| \leq 2^{- k}$ and $\gamma \geq \alpha$, we bound $2^{- k \alpha} (2^{-k}+|h|)^{\gamma-\alpha} \leq 2^{- k \gamma}$.
	Thus:
	\begin{equation}
		|u_k''|\leq \|\psi\|_{L^{\tilde{p}}}\|\check{\varphi}\|_{L^\infty} 2^{-k\gamma} \int_{B \left( 0, \epsilon_k \right)}f_K(k,h)2^{kd}\,dh \leq \|\psi\|_{L^{\tilde{p}}}\|\check{\varphi}\|_{L^\infty} 2^{-k\gamma} m_{K}^{(1)}(k).
	\end{equation}
	Denoting $\tilde{q}$ to be the H\"older conjugate of $q$, we have by H\"older's inequality 
	\begin{equation}
		\sum_{k} |u_k''|\leq \|\psi\|_{L^{\tilde{p}}}\|\check{\varphi}\|_{L^\infty}\left\|\left(2^{-k\gamma}\right)_{k \in \mathbb{N}}^{}\right\|_{\ell^{\tilde{q}}}\left\|m_{K}^{(1)}\right\|_{\ell^q} , 
	\end{equation}
\noindent which is finite because $\gamma > 0$ and $m_K^{\left( 1 \right)} \in \ell^q$.

Finally, \ref{item:convergence_3} follows immediately from the established estimates.
\end{proof}

\subsection{Reconstruction bound for \texorpdfstring{$\gamma> 0$}{gamma > 0}}

In \autoref{prop:convergence_approximating_distributions} just above, we have established the existence of a distribution $\mathcal{R} ( F ) \in \mathcal{D}^{\prime} ( \bar{K}_1 )$ that is a natural candidate for the reconstruction of the germ $F$.
In this section we focus on establishing that $\mathcal{R} ( F )$ does indeed satisfy the following reconstruction bound (recall the discussion of \autoref{def:Besov_reconstruction} and \autoref{Prop:Uniqueness}):
		\begin{equation}
			\left\| \left\| \sup\limits_{\psi \in \mathscr{B}^{r}} \left| \frac{\left( \mathcal{R} \left( F \right) - F_x \right) \left( \psi_x^{2^{-n}} \right)}{2^{- n \gamma}} \right| \right\|_{L^{p} \left( x \in K \right)} \right\|_{\ell^q \left( n \right)} < + \infty .
			\label{eq:reconstruction_bound_gamma>0}
		\end{equation}

In fact, we shall show that the left-hand term of \eqref{eq:reconstruction_bound_gamma>0} is bounded by a constant times $\left\| F \right\|_{\mathcal{G}_{p, q, q_1, \bar{K}_2, \varphi}^{\alpha, \beta, \gamma}}$.
For simplicity of notation, we denote: 
\begin{displaymath}
	G_{w}(\psi):=(\mathcal{R}(F)-F_{w})(\psi)\quad \text{for $w\in \R^d$ and $\psi \in \mathcal{D}(\bar{K}_1)$}. 
\end{displaymath}

As in the previous section, we shall be able to discuss $G$ thanks to the dyadic decomposition provided by the mollifier \eqref{EQN:MOLLIFIER}.
Remember that we defined $\mathcal{R} (F )$ as:
	\begin{equation}
		\mathcal{R} (F ) \left( \psi \right) \coloneqq \int_{\R^d} F_{z}\left (\hat{\varphi}^{\epsilon_{n_0}}*\hat{\varphi}_z^{2\epsilon_{n_0}}\right )\psi(z)\,dz + \sum\limits_{k = n_0}^{+ \infty} \int_{\R^d}F_{z}\left (\hat{\varphi}^{\epsilon_k}*\check{\varphi}_z^{\epsilon_k}\right )\psi(z)\,dz ,
	\end{equation}

\noindent where the right-hand term does not depend on the choice of $n_0 \in \mathbb{Z}$. Now for any fixed $w \in \mathbb{R}^d$, it holds by classical mollification of the distribution $F_w$ (with the mollifier $\rho$ defined as \eqref{EQN:MOLLIFIER}) that for any $n_0 \in \mathbb{Z}$:
	\begin{equation}
		F_w \left( \psi \right) = \int_{\R^d} F_{w}\left (\hat{\varphi}^{\epsilon_{n_0}}*\hat{\varphi}_z^{2\epsilon_{n_0}}\right )\psi(z)\,dz + \sum\limits_{k = n_0}^{+ \infty} \int_{\R^d}F_{w}\left (\hat{\varphi}^{\epsilon_k}*\check{\varphi}_z^{\epsilon_k}\right )\psi(z)\,dz .
	\end{equation}

Thus, we can decompose $G$ as:
	\begin{align}
		G_{w}(\psi) &= \int_{\R^d} \left( F_{z} - F_w \right) \left (\hat{\varphi}^{\epsilon_{n_0}}*\hat{\varphi}_z^{2\epsilon_{n_0}}\right )\psi(z)\,dz \\
		& \quad + \sum\limits_{k = n_0}^{+ \infty} \int_{\R^d} \left( F_{z} - F_w \right) \left (\hat{\varphi}^{\epsilon_k}*\check{\varphi}_z^{\epsilon_k}\right )\psi(z)\,dz
	\end{align}

Recalling \eqref{eq:distribution_against_convolution}, writing $F_z-F_w=(F_z-F_x)+(F_{x}-F_w)$ and taking into consideration the support of $\check{\varphi}_{z}^{\epsilon_n}$ and $\psi^{\lambda}_{w}$ (recall that $\supp ( \hat{\varphi} ) \subset B ( 0, 1 / 2 )$ and $\supp ( \check{\varphi} ) \subset B ( 0, 1 )$), we obtain, for $\psi \in \mathscr{B}^r$, $\lambda \in \left( 0, 1 \right]$ and any $n \in \mathbb{N}$, the more refined decomposition: 
\begin{align*}
	G_{w}(\psi^{\lambda}_{w})&=\int_{B(w,\lambda)}\int_{B(z,\epsilon_n)}(F_{z}-F_{x})\left (\hat{\varphi}^{\epsilon_n}_{x}\right )\hat{\varphi}^{\epsilon_{n - 1}}_{z}(x)\psi^{\lambda}_{w}(z)\,dx\,dz\\
	&+\int_{B(w,\lambda)}\int_{B(z,\epsilon_n)}(F_{x}-F_{w})\left (\hat{\varphi}^{\epsilon_n}_{x}\right )\hat{\varphi}^{\epsilon_{n - 1}}_{z}(x)\psi^{\lambda}_{w}(z)\,dx\,dz\\
	&+\sum_{k=n}^{\infty}\int_{B(w,\lambda)}\int_{B(z,\epsilon_k)}(F_{z}-F_{x})\left (\hat{\varphi}^{\epsilon_k}_{x}\right )\check{\varphi}^{\epsilon_k}_{z}(x)\psi^{\lambda}_{w}(z)\,dx\,dz\\
	&+\sum_{k=n}^{\infty}\int_{B(w,\lambda)}\int_{B(z,\epsilon_k)}(F_{x}-F_{w})\left (\hat{\varphi}^{\epsilon_k}_{x}\right )\check{\varphi}^{\epsilon_k}_{z}(x)\psi^{\lambda}_{w}(z)\,dx\,dz.
	\end{align*}
	Choosing $\lambda=2^{-n} = \epsilon_n$ gives us four terms which we define in the following way for $n\in \N$ and $w\in \R^d$:
	\begin{align}
		a_n(w)&:=\sup_{\psi \in \mathscr{B}^r}\left |\int_{B(w,\epsilon_n)}\int_{B(z,\epsilon_n)}(F_{z}-F_{x})\left (\hat{\varphi}^{\epsilon_n}_{x}\right )\hat{\varphi}^{\epsilon_{n - 1}}_{z}(x)\psi^{\epsilon_n}_{w}(z)\,dx\,dz\right |, \nonumber\\
		b_n(w)&:=\sup_{\psi \in \mathscr{B}^r}\left |\int_{B(w,\epsilon_n)}\int_{B(z,\epsilon_n)}(F_{x}-F_{w})\left (\hat{\varphi}^{\epsilon_n}_{x}\right )\hat{\varphi}^{\epsilon_{n - 1}}_{z}(x)\psi^{\epsilon_n}_{w}(z)\,dx\,dz\right |,\nonumber\\
		c_n(w)&:=\sup_{\psi \in \mathscr{B}^r}\left |\sum_{k=n}^{\infty}\int_{B(w,\epsilon_n)}\int_{B(z,\epsilon_k)}(F_{z}-F_{x})\left (\hat{\varphi}^{\epsilon_k}_{x}\right )\check{\varphi}^{\epsilon_k}_{z}(x)\psi^{\epsilon_n}_{w}(z)\,dx\,dz\right |,\nonumber\\
		d_n(w)&:=\sup_{\psi \in \mathscr{B}^r}\left |\sum_{k=n}^{\infty}\int_{B(w,\epsilon_n)}\int_{B(z,\epsilon_k)}(F_{x}-F_{w})\left (\hat{\varphi}^{\epsilon_k}_{x}\right )\check{\varphi}^{\epsilon_k}_{z}(x)\psi^{\epsilon_n}_{w}(z)\,dx\,dz\right |,\nonumber
	\end{align}

Now, the reconstruction bound is established in the following proposition. 
	\begin{proposition}\label{PROP:GAMMA>0}
		Assume that $m_K^{\left( 1 \right)}$, $m_K^{\left( 2 \right)}$, $m_K^{\left( 3 \right)}$ $\in \ell^q$ and that $\gamma \geq \alpha$. Then
			\begin{align*}
				&\left (\frac{\|a_n(w)\|_{L^{p}(w \in K)}}{2^{-n\gamma}}\right )_{n\in \N},\left (\frac{\|b_n(w)\|_{L^{p}(w \in K)}}{2^{-n\gamma}}\right )_{n\in \N},\\
				&\left (\frac{\|c_n(w)\|_{L^{p}(w \in K)}}{2^{-n\gamma}}\right )_{n\in \N},\left (\frac{\|d_n(w)\|_{L^{p}(w \in k)}}{2^{-n\gamma}}\right )_{n\in \N}\in \ell^q.
			\end{align*}
			As an immediate consequence, \eqref{eq:reconstruction_bound_gamma>0} holds, i.e.\ $\mathcal{R} \left( F \right)$ is a $r$-uniform $\gamma, p, q$-reconstruction of $F$ on $K$. 
	\end{proposition}

	\begin{proof}
		To begin, we focus on $\left(\frac{\|a_n(w)\|_{L^{p}(w \in K)}}{2^{-n\gamma}}\right)_{n\in \N}$. 
		
		By definition of $a_n$,
			\begin{equation}
				a_n(w) \leq \sup_{\psi \in \mathscr{B}^r} \int_{B(w,\epsilon_n)}\int_{B(z,\epsilon_n)} \left| (F_{z}-F_{x})\left (\hat{\varphi}^{\epsilon_n}_{x}\right ) \right| \left| \hat{\varphi}^{\epsilon_{n - 1}}_{z}(x) \right| \left| \psi^{\epsilon_n}_{w}(z) \right| \,dx\,dz .
			\end{equation}
			
		Using the estimates $| \hat{\varphi}^{\epsilon_{n - 1}}_{z}(x) | \leq 2^{n d} \| \hat{\varphi} \|_{L^{\infty}}$ and $\left| \psi^{\epsilon_n}_{w}(z) \right| \leq 2^{n d} \| \psi \|_{\mathcal{C}^r}$, we get:
			\begin{equation}
				a_n ( w ) \lesssim 2^{2 n d} \int_{B(w,\epsilon_n)}\int_{B(z,\epsilon_n)} \left| (F_{z}-F_{x})\left (\hat{\varphi}^{\epsilon_n}_{x}\right ) \right| \,dx\,dz .
			\end{equation}
			
		We apply the substitution $\tilde{x} = - ( x - z )$ then $\tilde{z} = z - w$ in this integral, which yield:
			\begin{equation}
				a_n ( w ) \lesssim 2^{2 n d} \int_{B(0,\epsilon_n)}\int_{B(0,\epsilon_n)} \left| (F_{\tilde{z} + w}-F_{- \tilde{x} + \tilde{z} + w})\left (\hat{\varphi}^{\epsilon_n}_{- \tilde{x} + \tilde{z} + w}\right ) \right| \,d \tilde{x} \,d \tilde{z} .
			\end{equation}		
			
		Now Minkowski's inequality implies:
			\begin{align}
				& \left\| a_n ( w ) \right\|_{L^p \left( w \in K \right)} \\ 
				& \quad \lesssim 2^{2 n d} \int_{B(0,\epsilon_n)}\int_{B(0,\epsilon_n)} \left\| (F_{\tilde{z} + w}-F_{- \tilde{x} + \tilde{z} + w})\left (\hat{\varphi}^{\epsilon_n}_{- \tilde{x} + \tilde{z} + w}\right ) \right\|_{L^p \left( w \in K \right)} \,d \tilde{x} \,d \tilde{z} .
			\end{align}					
		
		Applying the substitution $\tilde{w} = w + \tilde{z} - \tilde{x}$ in the $L^p$ norm yields:	
			\begin{align}
				\left\| a_n ( w ) \right\|_{L^p \left( w \in K \right)} & \lesssim 2^{2 n d} \int_{B(0,\epsilon_n)}\int_{B(0,\epsilon_n)} \left\| (F_{\tilde{w} + \tilde{x}}-F_{\tilde{w}})\left (\hat{\varphi}^{\epsilon_n}_{\tilde{w}}\right ) \right\|_{L^p \left( \tilde{w} \in \bar{K}_2 \right)} \,d \tilde{x} \,d \tilde{z} \\
				& = 2^{n d} \int_{B(0,\epsilon_n)} \left\| (F_{\tilde{w} + \tilde{x}}-F_{\tilde{w}})\left (\hat{\varphi}^{\epsilon_n}_{\tilde{w}}\right ) \right\|_{L^p \left( \tilde{w} \in \bar{K}_2 \right)} \,d \tilde{x} .
			\end{align}
		
		By definition of $f$, 
			\begin{equation}
				\left\| (F_{\tilde{w} + \tilde{x}}-F_{\tilde{w}})\left (\hat{\varphi}^{\epsilon_n}_{\tilde{w}}\right ) \right\|_{L^p \left( \tilde{w} \in \bar{K}_2 \right)} = f_K \left( n, \tilde{x} \right) 2^{- n \alpha} \left( 2^{-n} + \left| \tilde{x} \right| \right)^{\gamma - \alpha} .
			\end{equation}
		
		Since $\left| \tilde{x} \right| \leq 2^{-n}$ in the integral and $\gamma \geq \alpha$, this implies:
			\begin{equation}
				\left\| a_n ( w ) \right\|_{L^p \left( w \in K \right)} \lesssim 2^{n d} \int_{B(0,\epsilon_n)} f_K \left( n , \tilde{x} \right) 2^{- n \gamma} \,d \tilde{x} \leq 2^{- n \gamma} m_K^{\left( 1 \right)} \left( n \right) .
			\end{equation}						
			
		The assertion on $a_n$ follows.
		
		We now focus on $\left (\frac{\|b_n(w)\|_{L^{p}(w \in K)}}{2^{-n\gamma}}\right )_{n\in \N}$.
		By definition of $b_n$:
			\begin{equation}
				b_n \left( w \right) \leq \sup_{\psi \in \mathscr{B}^r} \int_{B(w,\epsilon_n)}\int_{B(z,\epsilon_n)} \left| (F_{x}-F_{w})\left (\hat{\varphi}^{\epsilon_n}_{x}\right ) \right| \left| \hat{\varphi}^{\epsilon_{n - 1}}_{z}(x) \right| \left| \psi^{\epsilon_n}_{w}(z) \right| \,dx\,dz .
			\end{equation}
			
		Once again, since $| \hat{\varphi}^{\epsilon_{n - 1}}_{z}(x) | \leq 2^{n d} \| \hat{\varphi} \|_{L^{\infty}}$ and $\left| \psi^{\epsilon_n}_{w}(z) \right| \leq 2^{n d} \| \psi \|_{\mathcal{C}^r}$, we obtain:
			\begin{equation}
				b_n \left( w \right) \lesssim 2^{2 n d} \int_{B(w,\epsilon_n)}\int_{B(z,\epsilon_n)} \left| (F_{x}-F_{w})\left (\hat{\varphi}^{\epsilon_n}_{x}\right ) \right| \,dx\,dz .
			\end{equation}
			
		Observe that for every $z\in B(w,\epsilon_n)$, we have $B(z,\epsilon_n) \subset B \left( w, 2 \epsilon_n \right)$. In turn we have
			\begin{equation}
				b_n \left( w \right) \lesssim 2^{n d} \int_{B(w, \epsilon_{n - 1})} \left| (F_{x}-F_{w})\left (\hat{\varphi}^{\epsilon_n}_{x}\right ) \right| \,dx .
			\end{equation}
			
		Substituting $\tilde{x} = - \left( x - w \right)$ then applying Minkowski's inequality yields:
			\begin{equation}
				\left\| b_n \left( w \right) \right\|_{L^p \left( w \in K \right)} \lesssim 2^{nd}\int_{B(0, \epsilon_{n - 1})} \left\| (F_{w - \tilde{x}}-F_{w})\left (\hat{\varphi}^{\epsilon_n}_{w - \tilde{x}}\right ) \right\|_{L^p \left( w \in K \right)} \,d \tilde{x} .
			\end{equation}
		
		Substituting $\tilde{w} = w - \tilde{x}$ in the $L^p$ norm yields:
			\begin{equation}
				\left\| b_n \left( w \right) \right\|_{L^p \left( w \in K \right)} \lesssim 2^{nd}\int_{B(0, \epsilon_{n - 1})} \left\| (F_{\tilde{w} + \tilde{x}}-F_{\tilde{w}})\left (\hat{\varphi}^{\epsilon_n}_{\tilde{w}}\right ) \right\|_{L^p \left( \tilde{w} \in \bar{K}_2 \right)} \,d \tilde{x} .
			\end{equation}
		
		Recalling the definition of $f_K$ and using the facts that $\left| \tilde{x} \right| \leq \epsilon_{n - 1}$ and $\gamma \geq \alpha$:
		\begin{equation}
			\left\| b_n \left( w \right) \right\|_{L^p \left( w \in K \right)} \lesssim 2^{- n \gamma} \int_{B(0, \epsilon_{n - 1})} 2^{n d} f_K \left( n, \tilde{x} \right) \,d \tilde{x} = 2^{- n \gamma} m_K^{\left( 1 \right)} \left( n \right) .
		\end{equation}
		
		The assertion on $b_n$ follows.
		
		Let us now consider $\left(\frac{\|c_n(w)\|_{L^{p}(w \in K)}}{2^{-n\gamma}}\right)_{n\in \N}$. 
		By definition of $c_n$:
			\begin{equation}
				c_n(w) \leq \sup_{\psi \in \mathscr{B}^r} \sum_{k=n}^{\infty}\int_{B(w,\epsilon_n)}\int_{B(z,\epsilon_k)} \left| (F_{z}-F_{x})\left (\hat{\varphi}^{\epsilon_k}_{x}\right )\right| \left| \check{\varphi}^{\epsilon_k}_{z}(x) \right| \left| \psi^{\epsilon_n}_{w}(z) \right| \,dx\,dz .
			\end{equation}
		
		Once again, since $| \check{\varphi}^{\epsilon_{k}}_{z}(x) | \leq 2^{k d} \| \check{\varphi} \|_{L^{\infty}}$ and $\left| \psi^{\epsilon_n}_{w}(z) \right| \leq 2^{n d} \| \psi \|_{\mathcal{C}^r}$, we obtain:
			\begin{equation}
				c_n(w) \lesssim 2^{n d} \sum_{k=n}^{\infty} 2^{k d} \int_{B(w,\epsilon_n)}\int_{B(z,\epsilon_k)} \left| (F_{z}-F_{x})\left (\hat{\varphi}^{\epsilon_k}_{x}\right )\right| dx\,dz .
			\end{equation}
		
		Reasoning as for $a_n$ we obtain:
			\begin{equation}
				\left\| c_n(w) \right\|_{L^p \left( w \in K \right)} \lesssim \sum_{k=n}^{\infty} 2^{k d} \int_{B(0,\epsilon_k)} 2^{-k \gamma} f_K \left( k, \tilde{x} \right) d \tilde{x} .
			\end{equation}

		Thus, 
			\begin{equation}
				\left\| c_n(w) \right\|_{L^p \left( w \in K \right)} \lesssim 2^{- n \gamma} m_K^{\left( 2 \right)} \left( n \right) .
			\end{equation}					
		
		The assertion on $c_n$ follows.
		
		Let us now consider $\left(\frac{\|d_n(w)\|_{L^{p}(w)}}{2^{-n\gamma}}\right)_{n\in \N}$. 
		Applying the substitution $\tilde{z}=z-w$ then $\tilde{x}=- \left( x-w \right)$, and remarking that the obtained integrand is supported in $\tilde{z} \in B \left( 0, \epsilon_n \right)$, $\tilde{x} \in B \left( - \tilde{z}, \epsilon_k \right)$, we can integrate over $\tilde{z} \in B \left( 0, \epsilon_n \right)$, $\tilde{x} \in B \left( 0, \epsilon_{n - 1} \right)$ without changing the value of the integral so that:
		\begin{align*}
			& d_n(w) \nonumber \\
			&  =\sup_{\psi \in \mathscr{B}^r} \Bigg|\sum_{k=n}^{\infty}\int_{B(0,\epsilon_n)}\int_{B(0,\epsilon_{n - 1})}(F_{-\tilde{x}+w}-F_{w})\left (\hat{\varphi}^{\epsilon_{k}}_{-\tilde{x}+w}\right ) \check{\varphi}^{\epsilon_k}_{\tilde{z}}(-\tilde{x})\psi^{\epsilon_n}(\tilde{z})\,d\tilde{x}\,d\tilde{z}\Bigg| \nonumber \\
			&  = \sup_{\psi \in \mathscr{B}^r} \Bigg|\sum_{k=n}^{\infty} \int_{B(0,\epsilon_{n - 1})} \left( (F_{-\tilde{x}+w}-F_{w})\left (\hat{\varphi}^{\epsilon_{k}}_{-\tilde{x}+w}\right ) \right) \left( \check{\varphi}^{\epsilon_k} * \psi^{\epsilon_n}(-\tilde{x}) \right) \,d\tilde{x} \Bigg| \nonumber \\
			&  \leq \sup_{\psi \in \mathscr{B}^r} \sum_{k=n}^{\infty} \int_{B(0,\epsilon_{n - 1})} \left| (F_{-\tilde{x}+w}-F_{w})\left (\hat{\varphi}^{\epsilon_{k}}_{-\tilde{x}+w}\right ) \right| \left| \check{\varphi}^{\epsilon_k} * \psi^{\epsilon_n}(-\tilde{x}) \right| \,d\tilde{x} \nonumber .
		\end{align*}
	
	Now from \autoref{LEMMA:CONV_INEQ}, $ \left| \check{\varphi}^{\epsilon_k} * \psi^{\epsilon_n}(-\tilde{x}) \right| \leq 2^{n \left( r + d \right) - k r} \|\psi\|_{C^r}\|\check{\varphi}\|_{L^1}$ so that:
	 \begin{equation}
	 	d_n \left( w \right) \lesssim \sum_{k=n}^{\infty} 2^{n \left( r + d \right) - k r} \int_{B(0,\epsilon_{n - 1})} \left| (F_{-\tilde{x}+w}-F_{w})\left (\hat{\varphi}^{\epsilon_{k}}_{-\tilde{x}+w}\right ) \right| \,d\tilde{x} .
	 \end{equation}
	
	Taking the $L^{p} \left( w \right)$ norm, applying Minkowski's inequality and substituting $\tilde{w} = w - \tilde{x}$ yields:
		\begin{equation}
			\left\| d_n \left( w \right) \right\|_{L^p \left( w \in K \right)} \lesssim \sum_{k=n}^{\infty} 2^{n \left( r + d \right) - k r} \int_{B(0,\epsilon_{n - 1})} \left\| (F_{\tilde{w} + \tilde{x}}-F_{\tilde{w}})\left (\hat{\varphi}^{\epsilon_{k}}_{\tilde{w}}\right ) \right\|_{L^p \left( \tilde{w} \in \bar{K}_2 \right)} \,d\tilde{x} .
		\end{equation}

	By definition of $f_K$:
		\begin{align}
			\left\| d_n \left( w \right) \right\|_{L^p \left( w \in K \right)} & \lesssim \sum_{k=n}^{\infty} 2^{n \left( r + d \right) - k r} \int_{B(0,\epsilon_{n - 1})} f_K \left( k, \tilde{x} \right) 2^{- k \alpha} 2^{- n \left( \gamma - \alpha \right)} \,d\tilde{x} \\
			& = 2^{- n \gamma} m_K^{\left( 3 \right)} \left( n \right) .
		\end{align}
	This is enough to conclude.
	\end{proof}
\begin{remark}\label{remark:d_n}
	Note that we did not use the assumption $\gamma > 0$ in these calculations.
\end{remark}

\subsection{The reconstruction for \texorpdfstring{$\gamma\leq 0$}{gamma <= 0}}

Now that we have treated the case $\gamma > 0$, let us discuss the problem of reconstruction when $\gamma \leq 0$.
It is very natural a priori to consider the same sequence of approximating distributions as in \autoref{DEF:APPROX}.
However, note from  \autoref{prop:convergence_approximating_distributions}, \autoref{item:convergence_2}, that the convergence of those approximating distributions fundamentally requires $\gamma > 0$. Namely, in this case, we cannot control the series $\sum_{k} u_k^{\prime \prime}$.

The idea, as in \cite{caravenna2020hairers}, is to simply remove the term $u_k^{\prime \prime}$ from the approximating sequence of \autoref{DEF:APPROX}. 
In particular, from \autoref{prop:convergence_approximating_distributions} we can still define:
	\begin{equation}
	\mathcal{R}(F)(\psi):=\int_{\R^d} F_{z}\left (\hat{\varphi}^{}*\hat{\varphi}_z^{2}\right )\psi(z)\,dz+\sum_{k=0}^{\infty}u'_k ,
	\label{eq:definition_reconstruction_gamma_negative}
	\end{equation}
	where we recall that $u_k^{\prime}$ is defined as:
	\begin{displaymath}
	u'_k:=\int_{\R^d}F_{x}\left(\hat{\varphi}^{\epsilon_k}_{x}\right)\left (\check{\varphi}^{\epsilon_k}*\psi\right )(x)\,dx, \quad \text{for $k\in \N$}.
	\end{displaymath}

Note that without the term $u_k^{\prime \prime}$, there is no simplification allowing to start the decomposition at a scale $2^{-n_0}$ for any $n_0$. 
Hence we need to take into account the fact that the sum starts at index $k = 0$ in \eqref{eq:definition_reconstruction_gamma_negative}.	
	
It remains to establish the reconstruction bound. As in the previous section, define for $w\in \R^d$ and $\psi \in \mathcal{D}(\R^d)$, $G_{w}(\psi) \coloneqq (\mathcal{R}(F)-F_{w})(\psi)$. 
Then, in a similar way to the previous section, it is straightforward to obtain the following decomposition:
	\begin{align*}
	G_w \left( \psi_w^{\epsilon_n} \right) & = \int_{B \left( w, \epsilon_n \right)} \int_{B \left( z, 1 \right)} \left( F_z - F_x \right) \left( \hat{\varphi}_x \right) \hat{\varphi}^2_z \left( x \right) \psi_w^{\epsilon_n} \left( z \right)\, d x\, d z \\
	& \quad +\int_{B \left( w, \epsilon_n \right)} \int_{B \left( z, 1 \right)} \left( F_x - F_w \right) \left( \hat{\varphi}_x \right) \hat{\varphi}^2_z \left( x \right) \psi_w^{\epsilon_n} \left( z \right)\, d x\, d z \\
	& \quad + \sum\limits_{k = 0}^{n - 1} \int_{B \left( w, \epsilon_n \right)} \int_{B \left( z, \epsilon_k \right)} \left( F_x - F_w \right) \left( \hat{\varphi}_x^{\epsilon_k} \right) \check{\varphi}_z^{\epsilon_k} \left( x \right) \psi_w^{\epsilon_n} \left( z \right)\, d x \,d z \\
	& \quad + \sum\limits_{k = n}^{+ \infty} \int_{B \left( w, \epsilon_n \right)} \int_{B \left( z, \epsilon_k \right)} \left( F_x - F_w \right) \left( \hat{\varphi}_x^{\epsilon_k} \right) \check{\varphi}_z^{\epsilon_k} \left( x \right) \psi_w^{\epsilon_n} \left( z \right) d x d z ,
	\end{align*}
	
\noindent so that checking the reconstruction bound follows from estimating the following quantities:
	\begin{align*}
	a_n(w) & := \sup_{\psi\in \mathscr{B}^r}\Bigg |\int_{B \left( w, \epsilon_n \right)} \int_{B \left( z, 1 \right)} \left( F_z - F_x \right) \left( \hat{\varphi}_x \right) \hat{\varphi}^2_z \left( x \right) \psi_w^{\epsilon_n} \left( z \right)\, d x\, d z \Bigg |,\\
	b_n(w) & := \sup_{\psi\in \mathscr{B}^r}\Bigg |\int_{B \left( w, \epsilon_n \right)} \int_{B \left( z, 1 \right)} \left( F_x - F_w \right) \left( \hat{\varphi}_x \right) \hat{\varphi}^2_z \left( x \right) \psi_w^{\epsilon_n} \left( z \right)\, d x\, d z \Bigg |,\\
	c_n(w) & := \sup_{\psi\in \mathscr{B}^r}\Bigg |\sum\limits_{k = 0}^{n - 1} \int_{B \left( w, \epsilon_n \right)} \int_{B \left( z, \epsilon_k \right)} \left( F_x - F_w \right) \left( \hat{\varphi}_x^{\epsilon_k} \right) \check{\varphi}_z^{\epsilon_k} \left( x \right) \psi_w^{\epsilon_n} \left( z \right)\, d x \,d z \Bigg |,\\
	d_n(w)& :=\sup_{\psi\in \mathscr{B}^r}\Bigg |\sum\limits_{k = n}^{+ \infty} \int_{B \left( w, \epsilon_n \right)} \int_{B \left( z, \epsilon_k \right)} \left( F_x - F_w \right) \left( \hat{\varphi}_x^{\epsilon_k} \right) \check{\varphi}_z^{\epsilon_k} \left( x \right) \psi_w^{\epsilon_n} \left( z \right) \,d x\, d z \Bigg |.
	\end{align*}
	
Note that these quantities are different from those named with the same letters in the section corresponding to $\gamma > 0$.
The reconstruction bound is established in the following proposition:
\begin{proposition}
In the setting of this section, assume that $\tilde{m}_K^{\left( 3 \right)}$, $\tilde{m}_K^{\left( 4 \right)} \in \ell^q$ and that $\gamma \geq \alpha$, $\gamma \leq 0$. Let $k$ be the scaling function defined in \eqref{eq:definition_scaling_function}. Then:
	\begin{align*}
	&\left (\frac{\|a_n(w)\|_{L^{p}(w \in K)}}{k \left( 2^{-n} \right)}\right )_{n\in \N},\left (\frac{\|b_n(w)\|_{L^{p}(w \in K)}}{k \left( 2^{-n} \right)}\right )_{n\in \N},\\
	&\left (\frac{\|c_n(w)\|_{L^{p}(w \in K)}}{k \left( 2^{-n} \right)}\right )_{n\in \N},\left (\frac{\|d_n(w)\|_{L^{p}(w \in K)}}{k \left( 2^{-n} \right)}\right )_{n\in \N}\in \ell^q.
	\end{align*}	
	As an immediate consequence, $\mathcal{R} \left( F \right)$ is a $r$-uniform $\gamma, p, q$-reconstruction of $F$ on $K$. 
\end{proposition}

\begin{remark}
The proof below actually works for any scaling function $k$ such that $\left( \frac{1}{k \left( 2^{- n} \right)} \right)_{n \in \mathbb{N}} \in \ell^q$.
\end{remark}

\begin{proof}
	The proof follows similarly as in \autoref{PROP:GAMMA>0}. 
	The same calculations as in \autoref{PROP:GAMMA>0} allow to establish:
		\begin{equation}
			\begin{dcases}
				\|a_n(w)\|_{L^p(w \in K)} &\lesssim m_K^{\left( 1 \right)}(0) , \\
				\|b_n(w)\|_{L^p(w \in K)} &\lesssim m_K^{\left( 1 \right)}(0) , \\
				\|d_n(w)\|_{L^p(w \in K)} &\lesssim 2^{- n \gamma} m_K^{\left( 3 \right)}(n) .
			\end{dcases}
		\end{equation}
	
	Let us now focus on $c_n$.
	By definition of $c_n$:
		\begin{equation}
			c_n(w) \leq \sup_{\psi\in \mathscr{B}^r} \sum\limits_{k = 0}^{n - 1} \int_{B \left( w, \epsilon_n \right)} \int_{B \left( z, \epsilon_k \right)} \left| \left( F_x - F_w \right) \left( \hat{\varphi}_x^{\epsilon_k} \right) \right| \left| \check{\varphi}_z^{\epsilon_k} \left( x \right) \right| \left| \psi_w^{\epsilon_n} \left( z \right) \right| \, d x \,d z .
		\end{equation}
		
	Using the estimates $| \check{\varphi}^{\epsilon_{k}}_{z}(x) | \leq 2^{k d} \| \check{\varphi} \|_{L^{\infty}}$ and $\left| \psi^{\epsilon_n}_{w}(z) \right| \leq 2^{n d} \| \psi \|_{\mathcal{C}^r}$, we get:
		\begin{equation}
			c_n(w) \lesssim \sum\limits_{k = 0}^{n - 1} 2^{n d} 2^{k d} \int_{B \left( w, 2^{-n} \right)} \int_{B \left( z, 2^{-k} \right)} \left| \left( F_x - F_w \right) \left( \hat{\varphi}_x^{2^{-k}} \right) \right| \, d x \,d z .
		\end{equation}
		
	In this integral, $B \left( z, 2^{-k} \right) \subset B \left( w, 2^{- k + 1} \right)$ so that:
		\begin{equation}
			c_n(w) \lesssim \sum\limits_{k = 0}^{n - 1} 2^{k d} \int_{B \left( w, 2^{-k+1} \right)} \left| \left( F_x - F_w \right) \left( \hat{\varphi}_x^{2^{-k}} \right) \right| \, d x .
		\end{equation}
	
	Substituting $\tilde{x} = - \left( x - w \right)$ in this integral:
		\begin{equation}
			c_n(w) \lesssim \sum\limits_{k = 0}^{n - 1} 2^{k d} \int_{B \left( 0, 2^{-k+1} \right)} \left| \left( F_{w - \tilde{x}} - F_w \right) \left( \hat{\varphi}_{w - \tilde{x}}^{2^{-k}} \right) \right| \, d \tilde{x} .
		\end{equation}

	Now we take the $L^p$ norm in $w$, apply Minkowski's inequality and change variable $\tilde{w} = w - \tilde{x}$ in the $L^p$	norm:
		\begin{equation}
		 \left\| c_n(w) \right\|_{L^p \left( w \in K \right)} \lesssim \sum\limits_{k = 0}^{n - 1} 2^{k d} \int_{B \left( 0, 2^{-k+1} \right)} \left\| \left( F_{\tilde{w} + \tilde{x}} - F_{\tilde{w}} \right) \left( \hat{\varphi}_{\tilde{w}}^{2^{-k}} \right) \right\|_{L^{p} \left( \tilde{w} \in \bar{K}_2 \right)} \, d \tilde{x} .
		\end{equation}
	
	Finally, remembering the definition of $f_K$ and using the fact that $\left| \tilde{x} \right| \leq 2^{- k + 1}$ in this integral: 
		\begin{equation}
			\left\| c_n(w) \right\|_{L^p \left( w \in K \right)} \lesssim 2^{- n \gamma} m_K^{\left( 4 \right)} \left( n \right) .
		\end{equation}
		
	This is enough to conclude, even for the case $\gamma = 0$.
	Indeed, for the terms $a$ and $b$, we use the fact that $\left( \frac{1}{k \left( 2^{- n} \right)} \right)_{n \in \mathbb{N}} \in \ell^q$.
	And for the terms $c$ and $d$, we have just established that:
		\begin{equation}
			\begin{dcases}
			\frac{\|c_n(w)\|_{L^p(w \in K)}}{k \left( 2^{- n} \right)} &\lesssim \tilde{m}_K^{\left( 4 \right)}(n) , \\
			\frac{\|d_n(w)\|_{L^p(w \in K)}}{k \left( 2^{- n} \right)} &\lesssim \tilde{m}_K^{\left( 3 \right)}(n) .
			\end{dcases}
		\end{equation}
\end{proof}

\subsection{The reconstruction is Besov}

We now show that $\mathcal{R}(F)$ lies in a suitable Besov space. 
\begin{proposition}\label{prop:reconstruction_is_besov}
In the setting of this section:
\begin{enumerate}
		\item (Global version) Assume that $\left\| F \right\|_{\mathcal{G}_{p, q, q_1, \mathbb{R}^d, \varphi}^{\alpha, \beta, \gamma}} < + \infty$. Then:
		\begin{enumerate}[ref =\emph{(\alph*)}]
			\item\label{item:reconstruction_is_besov_1} If $\beta \wedge \gamma > 0$, then $\mathcal{R} ( F ) = 0$.
			\item\label{item:reconstruction_is_besov_2} If $\beta \wedge \gamma = 0$, then for all $\kappa > 0$, $\mathcal{R} \left( F \right) \in \mathcal{B}_{p, 1}^{- \kappa}$.
			\item\label{item:reconstruction_is_besov_3} If $\beta \wedge \gamma < 0$, then $\mathcal{R} \left( F \right) \in \mathcal{B}_{p, q_1 \vee q}^{\beta \wedge \gamma}$.
		\end{enumerate}
	
	\item (Local version) Assume that for all $K \subset \mathbb{R}^d$, $\left\| F \right\|_{\mathcal{G}_{p, q, q_1, K, \varphi}^{\alpha, \beta, \gamma}} < + \infty$.
	Then there exists a global distribution $\mathcal{R} \left( F \right) \in \mathcal{D}^{\prime} ( \mathbb{R}^d )$ satisfying \eqref{eq:definition_reconstruction_bound_in_n} over all $K \subset \mathbb{R}^d$ and:
		\begin{enumerate}
			\item If $\beta \wedge \gamma > 0$, then $\mathcal{R} ( F ) = 0$.
			\item If $\beta \wedge \gamma = 0$, then for all $\kappa > 0$, $\mathcal{R} \left( F \right) \in \mathcal{B}_{p, 1, \mathrm{loc}}^{- \kappa}$.
			\item If $\beta \wedge \gamma < 0$, then $\mathcal{R} \left( F \right) \in \mathcal{B}_{p, q_1 \vee q, \mathrm{loc}}^{\beta \wedge \gamma}$.
		\end{enumerate}
	\end{enumerate}

Furthermore, the reconstruction map is continuous in the following sense: let $\mathcal{B}$ denote $\mathcal{B}_{p, 1, K}^{- \kappa}$ if $\beta \wedge \gamma = 0$ and $\mathcal{B}_{p, q_1 \vee q, K}^{\beta \wedge \gamma}$ if $\beta \wedge \gamma < 0$, then:
	\begin{equation}\label{eq:reconstruction_is_besov_continuity_estimate}
		\left\| \mathcal{R} \left( F \right) \right\|_{\mathcal{B}} \lesssim \left\| F \right\|_{\mathcal{G}_{p, q, q_1, \bar{K}_2, \varphi}^{\alpha, \beta, \gamma}} .
	\end{equation}	
\end{proposition}

\begin{proof}
Let us first prove the global version of the result.
The item \ref{item:reconstruction_is_besov_1} follows from \autoref{Prop:Uniqueness}, \autoref{item:scaling_prop_3}.

Now we turn to \ref{item:reconstruction_is_besov_3}.
Recall, from the equivalent definition of a Besov space \autoref{prop:Besov_equivalent_norms}, \autoref{item:prop_equivalent_besov_norm_1}, that it is sufficient to show, denoting $r \coloneqq q_1 \vee q$:
	\begin{displaymath}
	\left\| \left\| \frac{\mathcal{R} \left( F \right) \left( \hat{\varphi}_x^{2^{- n}} \right)}{2^{- n \left( \beta \wedge \gamma \right)}} \right\|_{L^p \left( x \in \mathbb{R}^d \right)} \right\|_{\ell^{r}(n)} < + \infty .
	\end{displaymath}
	
	By the assumption and the reconstruction bound obtained in the previous sections, we know that:
	\begin{displaymath}
	\left\| \left\| \frac{F_x \left( \hat{\varphi}_x^{\epsilon_n} \right)}{2^{-n \beta}} \right\|_{L^p \left( x \right)} \right\|_{\ell^{q_1}(n)} + \left\| \left\| \frac{ \left( \mathcal{R} \left( F \right) - F_x \right) \left( \hat{\varphi}_x^{\epsilon_n} \right)}{2^{-n \gamma}} \right\|_{L^p \left( x \right)} \right\|_{\ell^q (n)} \lesssim \left\| F \right\|_{\mathcal{G}_{p, q, q_1, \mathbb{R}^d, \varphi}^{\alpha, \beta, \gamma}} .
	\end{displaymath}
	Using the fact that $\beta, \gamma \geq \beta \wedge \gamma$ and the embeddings $\ell^{q_1} \subset \ell^r$, $\ell^{q} \subset \ell^r$:
	\begin{displaymath}
	\left\| \left\| \frac{F_x \left( \hat{\varphi}_x^{\epsilon_n} \right)}{2^{-n \left( \beta \wedge \gamma \right)}} \right\|_{L^p \left( x \right)} \right\|_{\ell^{r}(n)} + \left\| \left\| \frac{ \left( \mathcal{R} \left( F \right) - F_x \right) \left( \hat{\varphi}_x^{\epsilon_n} \right)}{2^{-n \left( \beta \wedge \gamma \right)}} \right\|_{L^p \left( x \right)} \right\|_{\ell^{r} (n)} \lesssim \left\| F \right\|_{\mathcal{G}_{p, q, q_1, \mathbb{R}^d, \varphi}^{\alpha, \beta, \gamma}}.
	\end{displaymath}
The triangle inequality yields the announced result.

Finally, \ref{item:reconstruction_is_besov_2} immediately follows from \ref{item:reconstruction_is_besov_3} after noticing that the condition of homogeneity $g \in \ell^{q_1}$ for some $\beta > 0$ implies $g \in \ell^{1}$ for any $\beta^{\prime} < \beta$; that the reconstruction bound for $\gamma$ implies the reconstruction bound for any $\gamma^{\prime} < \gamma$, see \autoref{Prop:Uniqueness}, \autoref{item:scaling_prop_1}; and that $\mathcal{B}_{p, \infty}^{- \kappa / 2} \subset \mathcal{B}_{p, 1}^{- \kappa}$ for any $\kappa > 0$.

Now let us discuss the local version of the result.
A global reconstruction $\mathcal{R} ( F )$ can be built by localization, as in \cite[Section~11]{caravenna2020hairers}.
Then, properties \ref{item:reconstruction_is_besov_1}, \ref{item:reconstruction_is_besov_2}, \ref{item:reconstruction_is_besov_3}, and \eqref{eq:reconstruction_is_besov_continuity_estimate} follow from the same arguments as in the global case, using the local version of \autoref{prop:Besov_equivalent_norms}, \autoref{item:prop_equivalent_besov_norm_1}.
\end{proof}

This concludes the proof of \autoref{Prop:Existence_gamma_gq_0}.


\section{Proof of \autoref{corollary:reconstruction_with_coherence_homogeneity}, \autoref{thm:properties_of_the_reconstruction_map} from \autoref{Prop:Existence_gamma_gq_0}}
\label{section:proof_of_corollary}

In this section, we prove \autoref{corollary:reconstruction_with_coherence_homogeneity} and \autoref{thm:properties_of_the_reconstruction_map} from \autoref{Prop:Existence_gamma_gq_0}.

The case $\gamma \leq 0$ is treated similarly to the case $\gamma > 0$ so we only treat the case $\gamma > 0$ for concision.

Thus, we consider a germ $F$, and reals $p, q \in \left[ 1, + \infty \right]$, $\alpha, \beta, \gamma \in \mathbb{R}$ with $\alpha \leq \gamma$ and $\gamma > 0$.
We let $K \subset \mathbb{R}^d$ and we assume that there is a test-function $\varphi \in \mathcal{D} ( \mathbb{R}^d )$ with $\int \varphi \neq 0$ such that \eqref{eq:homogeneity_condition} and \eqref{eq:coherence_condition_1} hold.
Let $r \in \mathbb{N}$ be an integer such that $r > \max \left( - \alpha, - \beta \right)$.

We shall show that the hypotheses of \autoref{Prop:Existence_gamma_gq_0} are satisfied.
First, notice that this requires us to exhibit a test-function $\hat{\varphi} \in \mathcal{D} ( B (0, 1 / 2) )$ such that $\int \hat{\varphi} = 1$ and $\int x^k \hat{\varphi} \left( x \right) dx = 0$ for $1 \leq \left| k \right| \leq r - 1$.

For this purpose, we \emph{tweak} the test-function $\varphi$ as presented in \cite{caravenna2020hairers}.

\begin{lemma}\label{LEMMA:TWEAK}[Tweaking a test-function, \cite[Lemma 8.1]{caravenna2020hairers}]
	Fix $r\in \N$, distinct $\lambda_0,\lambda_1,...,\lambda_{r-1}\in(0,\infty)$ and a test-function $\varphi \in \mathcal{D}(\R^d)$ such that $\int \varphi \ne 0$ and $\supp{\varphi}\subset B(0,R_{\varphi})$. 

	Then there exists constants $c_0,c_1,....,c_{r-1}\in \R$ such that the \emph{tweaked} test-function $\hat{\varphi}$, defined by 
	\begin{equation}
	\hat{\varphi}:=\frac{1}{\int \varphi}\sum_{i=0}^{r-1}c_i\varphi^{\lambda_{i}},
	\end{equation}
	has the following properties:
	\begin{enumerate}
		\item\label{item:tweaking_1} $\int_{\R^d} \hat{\varphi} = 1$,
		\item\label{item:tweaking_2} $\hat{\varphi}$ annihilates monomials of degree 1 to $r-1$, specifically 
		\begin{displaymath}
		\int_{\R^d}y^k\hat{\varphi}(y)\,dy=0, \quad \text{for all $k\in\N_{0}^d:1\leq|k|\leq r-1$},
		\end{displaymath}
		\item moreover, if $0 < \lambda_{i} < \frac{1}{2 R_{\varphi}}$ for $i=0,1...,r-1$, then:
		\begin{displaymath}
		\supp \hat{\varphi}\subset B(0,1/2).
		\end{displaymath}
	\end{enumerate}
\end{lemma}
\begin{remark}
	Whilst we direct the reader to \cite[Lemma 8.1]{caravenna2020hairers} for the proof we outline the approach here. The main idea is to consider an arbitrary $\hat{\varphi}:=( 1/{\int \varphi} ) \sum_{i=0}^{r-1}c_i\varphi^{\lambda_{i}}$ and write a system of linear equations in $c_i$ from the relations given by \autoref{item:tweaking_1} and \autoref{item:tweaking_2} of \autoref{LEMMA:TWEAK}.
	The obtained system of equations involves the Vandermonde matrix of $(\lambda_i)_{0 \leq i \leq r - 1}$.
	Since the $\lambda_i$ are distinct, this system is invertible, and one can even provide explicit expressions for $c_i$ \cite[see~Equation~(8.1)]{caravenna2020hairers}.
\end{remark}

\begin{remark}
In \cite{caravenna2020hairers}, the authors choose the specific values $\lambda_i \coloneqq \frac{2^{- i - 1}}{1 + R_{\varphi}}$ for $i=0,1...,r-1$. 
The reason for this choice is that it allows for explicit quantitative bounds throughout their calculations.
In this paper however, we do not track the precise constants that appear in the estimates, which is why we do not pick explicit values for $\lambda_i$.
\end{remark}

Now we show that the properties of homogeneity and coherence are stable by tweaking.
\begin{proposition}
Assume that there exists $\varphi \in \mathcal{D} ( \mathbb{R}^d )$ such that $\int \varphi \neq 0$ and \eqref{eq:homogeneity_condition} resp.\ \eqref{eq:coherence_condition_1} is satisfied for $\varphi$.
Then, for any $r \in \mathbb{N}$, there exists $\hat{\varphi} \in \mathcal{D} ( B \left( 0, 1 / 2 \right))$ such that $\int \hat{\varphi} = 1$, $\int x^k \hat{\varphi} (x ) dx = 0$ for $1 \leq \left| k \right| \leq r - 1$ and \eqref{eq:homogeneity_condition} resp.\ \eqref{eq:coherence_condition_1} is satisfied for $\hat{\varphi}$.
\label{prop:stabilty_of_homogeneity_coherence_by_tweaking}
\end{proposition}

\begin{proof}
Let us only present the proof for the coherence condition \eqref{eq:coherence_condition_1} as the other case is similar to treat.
Thus, assume that $\varphi \in \mathcal{D} ( \mathbb{R}^d )$ is such that $\int \varphi \neq 0$ and:
	\begin{equation}
		 \left\| \left\| \left\| \frac{\left( F_{x + h} - F_x \right) \left( \varphi_x^{2^{-n}} \right)}{2^{- n \alpha} \left( 2^{-n} + \left| h \right| \right)^{\gamma - \alpha}} \right\|_{L^{p} \left( x \in \bar{K}_2 \right)} \right\|_{\ell^{\infty} \left( n \right)} \right\|_{L_h^q \left( h \in B \left( 0, 2 \right) \right)} < + \infty .
	\end{equation}
	
As in \autoref{LEMMA:TWEAK}, let $R_{\varphi} > 0$ be such that $\supp \left( \varphi \right) \subset B \left( 0, R_{\varphi} \right)$.
Let $r \in \mathbb{N}$ and fix distinct $m_0, \dots , m_{r - 1} \in \mathbb{N}$ such that for $0 \leq \left| k \right| \leq r - 1$, $2^{- m_k} < \frac{1}{2 R_{\varphi}}$.
We apply \autoref{LEMMA:TWEAK} to $\lambda_k \coloneqq 2^{-m_k}$. Let us denote 
	\begin{equation}
		\hat{\varphi}:=\frac{1}{\int \varphi}\sum_{i=0}^{r-1}c_i\varphi^{2^{- m_i}} ,
	\end{equation}
	
\noindent the obtained test-function, then $\hat{\varphi} \in \mathcal{D} ( B \left( 0, 1 / 2 \right))$ and $\int \hat{\varphi} = 1$, $\int x^k \hat{\varphi} (x ) dx = 0$ for $1 \leq \left| k \right| \leq r - 1$.
To conclude, it is enough by triangle inequality to show that for any $m \in \mathbb{N}$:
	\begin{equation}
		 \left\| \left\| \left\| \frac{\left( F_{x + h} - F_x \right) \left( \varphi_x^{2^{-n - m}} \right)}{2^{- n \alpha} \left( 2^{-n} + \left| h \right| \right)^{\gamma - \alpha}} \right\|_{L^{p} \left( x \in \bar{K}_2 \right)} \right\|_{\ell^{\infty} \left( n \right)} \right\|_{L_h^q \left( h \in B \left( 0, 2 \right) \right)} < + \infty .
	\end{equation}
	
But this follows immediately from the estimate:
	\begin{equation}
		2^{- n \alpha} \left( 2^{- n} + \left| h \right| \right)^{\gamma - \alpha} \geq 2^{m \alpha} 2^{- \left( n + m \right) \alpha} \left( 2^{- \left( n + m \right)} + \left| h \right| \right)^{\gamma - \alpha} .
	\end{equation}		
(Remember that $\gamma \geq \alpha$ from our hypotheses).
\end{proof}

We still have to translate the conditions of \autoref{corollary:reconstruction_with_coherence_homogeneity} into the conditions of \autoref{Prop:Existence_gamma_gq_0}.
This is possible thanks to the \autoref{LEMMA:tech}, which implies the following result.
Note that \autoref{LEMMA:tech} is a purely elementary result.
However, as its proof is a bit technical, we state it and prove it in the appendix. 

\begin{corollary}\label{corollary:translating_coherence_to_m_1_m_2_m_3}
Let $f \colon \mathbb{N} \times \mathbb{R}^d \to \mathbb{R}_+$ be a positive function, and $c > 0$. Define $m_f^{\left( 1 \right)}$, $m_{c, f}^{\left( 2 \right)}$, $m_{c, f}^{\left( 3 \right)}$ as in \eqref{eq:definition_of_m1_m2_m3}. Assume that:
			\begin{equation} \label{eq:technical_corollary_condition}
				\left\| \left\| f \left( k, x \right) \right\|_{\ell^\infty \left( k \right)} \right\|_{L_x^q \left( x \in B \left( 0, 2 \right) \right)} < + \infty .
			\end{equation}
Then $m_f^{\left( 1 \right)}$, $m_{c, f}^{\left( 2 \right)}$, $m_{c, f}^{\left( 3 \right)}$ are in $\ell^q$.	

\end{corollary}

\begin{proof}
Applying \autoref{LEMMA:tech} to $a_{k, n} \coloneqq \delta_{k, n}$ resp.\ $a_{k, n} \coloneqq 2^{- \left( k - n \right) c} \mathds{1}_{\left\lbrace k \geq n \right\rbrace}$ immediately gives the result for $m_f^{\left( 1 \right)}$ resp.\ $m_{c, f}^{\left( 2 \right)}$.

Now let us treat $m_{c, f}^{\left( 3 \right)}$. For $n \in \mathbb{N}$, $h \in \mathbb{R}^d$, set:
	\begin{equation}
		\tilde{f} \left( n, h \right) \coloneqq \sum\limits_{k = 0}^{+ \infty} 2^{- \left( k - n \right) c} \mathds{1}_{\left\lbrace k \geq n \right\rbrace} f \left( k, h \right) ,
	\end{equation}

\noindent so that $m_{c, f}^{\left( 3 \right)} = m_{\tilde{f}}^{\left( 1 \right)}$.
It is straightforward to see that if \eqref{eq:technical_corollary_condition} is satisfied for $f$, then \eqref{eq:technical_corollary_condition} is also satisfied for $\tilde{f}$.
We conclude by applying lemma \autoref{LEMMA:tech} to $a_{k, n} \coloneqq \delta_{k, n}$ and the function $\tilde{f}$.\end{proof}

Combining \autoref{prop:stabilty_of_homogeneity_coherence_by_tweaking}, \autoref{corollary:translating_coherence_to_m_1_m_2_m_3} and \autoref{Prop:Existence_gamma_gq_0} yields \autoref{corollary:reconstruction_with_coherence_homogeneity} and \autoref{thm:properties_of_the_reconstruction_map} (in the case $\gamma > 0$).

\section{Proof of \autoref{prop:Young_multiplication_besov} from \autoref{Prop:Existence_gamma_gq_0}}
\label{section:application_young_multiplication}

Let us now prove \autoref{prop:Young_multiplication_besov} from \autoref{Prop:Existence_gamma_gq_0}.
In the remainder of this section, we consider the setting of \autoref{subsection:young_multiplication}.
That is, we let $\alpha < 0, \beta > 0$ with $\alpha + \beta > 0$ and $\beta \notin \mathbb{N}$, as well as $p_1, p_2, q_1, q_2, p, q \in \left[ 1, + \infty \right]$ with $\frac{1}{p} = \frac{1}{p_1} + \frac{1}{p_2}$, $\frac{1}{q} = \frac{1}{q_1} + \frac{1}{q_2}$; we fix distributions $g \in \mathcal{B}_{p_1, q_1}^{\alpha}$, $f \in \mathcal{B}_{p_2, q_2}^{\beta}$; and we define germs $F$ and $P$ as in \eqref{eq:Taylor_germ_of_f} resp.\ \eqref{eq:germ_for_multiplication}.

Fix any test-function $\varphi \in \mathcal{D} ( \mathbb{R}^d )$.
Recalling the statement of \autoref{Prop:Existence_gamma_gq_0}, \autoref{prop:Young_multiplication_besov} holds as soon as the following quantities are finite (recall that $\epsilon_n \coloneqq 2^{-n}$):
	\begin{equation}
		\begin{dcases}
			v_1  \coloneqq \left\| \left\| \frac{P_x \left( \varphi_x^{\epsilon_n} \right)}{2^{- n \alpha}} \right\|_{L^{p} \left( x \right)} \right\|_{\ell^{q_1} \left( n \right)} , \\
			v_2 \coloneqq \left| \int_{h \in B \left( 0, \epsilon_{n - 1} \right)} 2^{n d} \left\| \frac{\left( P_{x + h} - P_x \right) \left( \varphi_x^{\epsilon_n} \right)}{2^{- n \alpha} \left( 2^{-n} + \left| h \right| \right)^{\beta}} \right\|_{L^{p} \left( x \right)} d h \right\|_{\ell^q \left( n \right)} , \\
			v_3 \coloneqq \left\| \sum\limits_{k = n}^{+ \infty} \int_{h \in B \left( 0, \epsilon_k \right)} 2^{- \left( k - n \right) \left( \alpha + \beta \right) + k d} \left\| \frac{\left( P_{x + h} - P_x \right) \left( \varphi_x^{\epsilon_k} \right)}{2^{- k \alpha} \left( 2^{-k} + \left| h \right| \right)^{\beta}} \right\|_{L^{p} \left( x \right)} d h \right\|_{\ell^q \left( n \right)} , \\
			v_4 \coloneqq \left\| \sum\limits_{k = n}^{+ \infty} \int_{h \in B \left( 0, \epsilon_{n - 1} \right)} 2^{- \left( k - n \right) \left( \alpha + r \right) + n d} \left\| \frac{\left( P_{x + h} - P_x \right) \left( \varphi_x^{\epsilon_k} \right)}{2^{- k \alpha} \left( 2^{-k} + \left| h \right| \right)^{\beta}} \right\|_{L^{p} \left( x \right)} d h \right\|_{\ell^q \left( n \right)} .
		\end{dcases}
		\label{eq:definition_of_g_m1_m2_m3}
	\end{equation}

This is the content of the following proposition.
\begin{proposition}
$v_1, v_2, v_3, v_4 < + \infty$.
As a consequence, \autoref{prop:Young_multiplication_besov} holds.
\end{proposition}

\begin{proof}
We start with $v_1$.
For $x \in \mathbb{R}^d, \varphi \in \mathcal{D} ( \mathbb{R}^d )$ we have:
	\begin{equation}
		P_x \left( \varphi \right) = \sum\limits_{0 \leq \left| k \right| < \beta} \frac{\partial^k f \left( x \right)}{k !} g \left( \left( \cdot - x \right)^k \varphi \left( \cdot \right) \right) .
	\end{equation}

Thus by triangle inequality and H\"older's inequality:
	\begin{align}
		\left\| \frac{P_x \left( \varphi_x^{\epsilon_n} \right)}{2^{- n \alpha}} \right\|_{L^p \left( x \right)} & \leq \sum\limits_{0 \leq \left| k \right| < \beta} \frac{\left\| \partial^k f \right\|_{L^{p_2}}}{k !} \frac{\left\| g \left( \left( \cdot - x \right)^k \varphi_x^{\epsilon_n} \left( \cdot \right) \right) \right\|_{L^{p_1} \left( x \right)}}{2^{- n \alpha}} .
	\end{align}

Now summing in $n$:
	\begin{equation}
		v_1 \leq \sum\limits_{0 \leq \left| k \right| < \beta} \frac{\left\| \partial^k f \right\|_{L^{p_2}}}{k !} \left\| \frac{\left\| g \left( \left( \cdot - x \right)^k \varphi_x^{\epsilon_n} \left( \cdot \right) \right) \right\|_{L^{p_1} \left( x \right)}}{2^{- \alpha n}} \right\|_{\ell^{q_1} \left( n \right)} .\nonumber
	\end{equation}		

Considering the test-function $\psi \colon z \mapsto z^k \phi \left( z \right)$ and recalling the definition of Besov spaces, \autoref{def:Besov_spaces_by_local_means}, one obtains:
	\begin{equation}
		\left\| \frac{\left\| g \left( \left( \cdot - x \right)^k \varphi_x^{\epsilon_n} \left( \cdot \right) \right) \right\|_{L^{p_1} \left( x \right)}}{2^{- \left( \alpha + | k | \right) n}} \right\|_{\ell^{q_1} \left( n \right)} \lesssim \left\| g \right\|_{\mathcal{B}_{p_1, q_1}^{\alpha}} .
	\end{equation}

Recall also that since $\beta > |k|$ in the sum above, one has $\left\| \partial_k f \right\|_{L^{p_2}} \lesssim \left\| \partial_k f \right\|_{\mathcal{B}_{p_2, q_2}^{\beta - |k |}} \lesssim \left\| f \right\|_{\mathcal{B}_{p_2, q_2}^{\beta}}$, and thus:
	\begin{equation}
		v_1 \lesssim \left\| g \right\|_{\mathcal{B}_{p_1, q_1}^{\alpha}} \left\| f \right\|_{\mathcal{B}_{p_2, q_2}^{\beta}} < + \infty .
	\end{equation}

Now we consider the quantities $v_2, v_3, v_4$.
Let us use the following notation for the Taylor expansions: if $f$ is a sufficiently regular function, $\alpha \in \mathbb{R}$ and $x, h \in \mathbb{R}^d$, set $T_f^{\alpha} \left( x, h \right) \coloneqq f \left( x + h \right) - \sum_{0 \leq \left| l \right| < \alpha} \frac{1}{l !} \partial^{l} f \left( x \right) h^l$. 

Let $x, h, z \in \mathbb{R}^d$ then a straightforward calculation establishes (recall that $F$ is defined in \eqref{eq:Taylor_germ_of_f}):
	\begin{align}
		\left( F_{x + h} - F_x \right) \left( z \right) & = -\sum\limits_{0 \leq \left| k \right| < \beta} \frac{\left( z - x \right)^k}{k !} T_{\partial^k f}^{\beta - \left| k \right|} \left( x, h \right) .
	\end{align}
	
For simplicity of notations, denote $T_k \left( x, h \right) \coloneqq T_{\partial^k f}^{\beta - \left| k \right|} \left( x, h \right)$ for the remainder of this proof.
Let $\varphi \in \mathcal{D} ( \mathbb{R}^d )$ be any test-function.
We deduce that for $x, h \in \mathbb{R}^d$:
	\begin{align}
		\left( P_{x + h} - P_x \right) \left( \varphi \right) & = g \left( \varphi \left( \cdot \right) \left( F_{x + h} - F_x \right) \left( \cdot \right) \right) \\
		& = -\sum\limits_{0 \leq \left| k \right| < \beta} \frac{1}{k !} T_k \left( x, h \right) g \left( \left( \cdot - x \right)^k \varphi \left( \cdot \right) \right) .
	\end{align}

Applying the triangle inequality then H\"older's inequality:
	\begin{align}
		& \left\| \left( P_{x + h} - P_x \right) \left( \varphi \right) \right\|_{L^p \left( x \right)} \leq \sum\limits_{0 \leq \left| k \right| < \beta} \frac{\left\| g \left( \left( \cdot - x \right)^k \varphi \left( \cdot \right) \right) \right\|_{L^{p_1} \left( x \right)} \left\| T_k \left( x, h \right) \right\|_{L^{p_2} \left( x \right)}}{k !} .
	\end{align}

For $0 \leq \left| k \right| < \beta$, it holds that $\lambda^{\alpha} \left( \lambda + \left| h \right| \right)^{\beta} \gtrsim \lambda^{\alpha + \left| k \right|} \left| h \right|^{\beta - \left| k \right|}$, so that:
	\begin{align}
		& \left\| \frac{\left( P_{x + h} - P_x \right) \left( \varphi_x^{\lambda} \right)}{\lambda^{\alpha} \left( \lambda + \left| h \right| \right)^{\beta}} \right\|_{L^p \left( x \right)} \\
		& \quad \lesssim \sum\limits_{0 \leq \left| k \right| < \beta} \frac{1}{k !} \frac{\left\| g \left( \left( \cdot - x \right)^k \varphi_x^{\lambda} \left( \cdot \right) \right) \right\|_{L^{p_1} \left( x \right)}}{\lambda^{\alpha + \left| k \right|}} \frac{\left\| T_k \left( x, h \right) \right\|_{L^{p_2} \left( x \right)}}{\left| h \right|^{\beta - \left| k \right|}} .
	\end{align}

For $0 \leq \left| k \right| < \beta$, set :
	\begin{equation}
		\begin{dcases}
			\mu_k \left( n \right) & \coloneqq \frac{1}{k !} \frac{\left\| g \left( \left( \cdot - x \right)^k \varphi_x^{\lambda} \left( \cdot \right) \right) \right\|_{L^{p_1} \left( x \right)}}{\lambda^{\alpha + \left| k \right|}} , \\
			\nu_k \left( h \right) & \coloneqq \frac{\left\| T_k \left( x, h \right) \right\|_{L^{p_2} \left( x \right)}}{\left| h \right|^{\beta - \left| k \right|}}
		\end{dcases}
	\end{equation}

Reasoning as above, it holds that $\mu_k \in \ell^{q_1} \left( n \in \mathbb{N} \right)$.
Also, recalling \autoref{prop:Besov_equivalent_norms}, one observes that $\nu_k \in L_h^{q_2} \left( h \in B ( 0, 2 ) \right)$.
We conclude by applying \autoref{LEMMA:tech}, \autoref{item:technical_lemma_condition_2} to the quantities $v_2, v_3, v_4$ in the same way as in the proof of \autoref{corollary:translating_coherence_to_m_1_m_2_m_3}.

We obtain \autoref{prop:Young_multiplication_besov} by setting $\mathcal{M} \left( g, f \right) \coloneqq \mathcal{R} \left( P \right)$.
Note that by collecting all the inequalities, we even obtain the following continuity estimate:
	\begin{equation}
		\left\| \mathcal{M} \left( g, f \right) \right\|_{\mathcal{B}_{p, q_1}^{\alpha}} \lesssim \left\| g \right\|_{\mathcal{B}_{p_1, q_1}^{\alpha}} \left\| f \right\|_{\mathcal{B}_{p_2, q_2}^{\beta}} .
	\end{equation}
\end{proof}


\appendix

\section{Besov spaces}
\label{appendix_besov_spaces}

In this section, we define Besov spaces, and recall some of their properties.
There are many different equivalent norms used in the literature to define and study Besov spaces.
In our context, the following definition \enquote{by local means} from \cite{MR3684891} will be the most useful.

\begin{definition}[Besov spaces]
Let $\alpha \in \mathbb{R}$, $p, q \in \left[ 1, + \infty \right]$. Let $r \in \mathbb{N}$ be such that $r > - \alpha$, and let $n_0 \in \mathbb{Z}$. We define $\mathcal{B}_{p, q}^{\alpha} = \mathcal{B}_{p, q}^{\alpha} ( \mathbb{R}^d )$ to be the space of distributions $f \in \mathcal{D}^{\prime} ( \mathbb{R}^d )$ such that:
	\begin{equation}
		\begin{dcases}
			\left\| \left\| \sup\limits_{\psi \in \mathscr{B}^{r}} \left| \frac{f \left( \psi_x^{2^{-n}} \right)}{2^{- n \alpha}} \right| \right\|_{L^p \left( x \right)} \right\|_{\ell^q \left( n \geq n_0 \right)} < + \infty & \text{if } \alpha < 0 , \\
			\left\| \sup\limits_{\psi \in \mathscr{B}^{r}} \left| f \left( \psi_x \right) \right| \right\|_{L^p \left( x \right)} + \left\| \left\| \sup\limits_{\psi \in \mathscr{B}_{\left\lfloor \alpha \right\rfloor}^{r}} \left| \frac{f \left( \psi_x^{2^{-n}} \right)}{2^{- n \alpha}} \right| \right\|_{L^p \left( x \right)} \right\|_{\ell^q \left( n \geq n_0 \right)} < + \infty & \text{if } \alpha \geq 0 .
		\end{dcases}
	\end{equation}

Here, recall that $\mathscr{B}_{\left\lfloor \alpha \right\rfloor}^{r}$ denotes the space of test-functions $\psi \in \mathscr{B}^{r}$ such that $\int x^k \psi \left( x \right) d x = 0$ for $0 \leq \left| k \right| \leq \left\lfloor \alpha \right\rfloor$.
\label{def:Besov_spaces_by_local_means}
\end{definition}

\begin{remark}\label{remark:does_not_depend_on_n0}
In \cite[Proposition~2.4]{MR3684891}, it is established (in the case $n_0 = 0$) that this definition does not depend on the choice of $r > - \alpha$ and that it is equivalent to the usual definition \enquote{by wavelets}.
It is also straightforward to establish that the definition does not depend on the choice of $n_0 \in \mathbb{Z}$, so that unless specified, $n_0$ is taken to be $0$.
\end{remark}

\begin{remark}\label{remark:Besov_spaces_tempered_distributions}
From \autoref{def:Besov_spaces_by_local_means}, we a priori only have $\mathcal{B}_{p, q}^{\alpha} \subset \mathcal{D}^{\prime}$ (the space of Schwartz distributions), while usual definitions of Besov spaces impose $\mathcal{B}_{p, q}^{\alpha} \subset \mathcal{S}^{\prime}$ (the space of tempered distributions).
However, the latter inclusion is actually a consequence of our definition, which can be seen for instance from the wavelet characterisation \cite[Proposition~2.4]{MR3684891}.
\end{remark}

In some situations, it is useful to have local versions of the spaces $\mathcal{B}_{p, q}^{\alpha} ( \mathbb{R}^d )$. In this case, the bounds of \autoref{def:Besov_spaces_by_local_means} are required to hold on $L^p ( x \in K )$ for all compact $K$, rather than on $L^p ( x \in \mathbb{R}^d )$.
\begin{definition}[Local Besov spaces]
Let $\alpha \in \mathbb{R}$, $p, q \in \left[ 1, + \infty \right]$. Let $r \in \mathbb{N}$ be such that $r > - \alpha$. We define $\mathcal{B}_{p, q, \mathrm{loc}}^{\alpha} = \mathcal{B}_{p, q, \mathrm{loc}}^{\alpha} ( \mathbb{R}^d )$ to be the space of distributions $f \in \mathcal{D}^{\prime} ( \mathbb{R}^d )$ such that for all compact $K \subset \mathbb{R}^d$:
	\begin{equation}
		\begin{dcases}
			\left\| \left\| \sup\limits_{\psi \in \mathscr{B}^{r}} \left| \frac{f \left( \psi_x^{2^{-n}} \right)}{2^{- n \alpha}} \right| \right\|_{L^p \left( x \in K \right)} \right\|_{\ell^q \left( n \right)} < \infty & \text{if } \alpha < 0 , \\
			\left\| \sup\limits_{\psi \in \mathscr{B}^{r}} \left| f \left( \psi_x \right) \right| \right\|_{L^p \left( x \in K \right)} + \left\| \left\| \sup\limits_{\psi \in \mathscr{B}_{\left\lfloor \alpha \right\rfloor}^{r}} \left| \frac{f \left( \psi_x^{2^{-n}} \right)}{2^{- n \alpha}} \right| \right\|_{L^p \left( x \in K \right)} \right\|_{\ell^q \left( n \right)} < \infty & \text{if } \alpha \geq 0 .
		\end{dcases}
	\end{equation}

We note $\left\| f \right\|_{\mathcal{B}_ {p, q, K}^{\alpha}}$ the norm provided by the quantity just above.	
Using the same argument as in \autoref{remark:does_not_depend_on_n0}, this definition does not depend on the choice of $r > - \alpha$.
\end{definition}

The following equivalent norms will be useful for us.

\begin{proposition}\label{prop:Besov_equivalent_norms} 
Let $\alpha \in \mathbb{R}$, $p, q \in \left[ 1, + \infty \right]$, and $f \in \mathcal{D}^{\prime} ( \mathbb{R}^d )$.
	\begin{enumerate}[ref =\emph{(\arabic*)}]
		\item\label{item:prop_equivalent_besov_norm_1} If $\alpha < 0$, then $f \in \mathcal{B}_{p, q}^{\alpha} ( \mathbb{R}^d )$ if and only if there exists a test-function $\varphi \in \mathcal{D} ( \mathbb{R}^d )$ such that $\int \varphi \neq 0$ and:
			\begin{equation}
				\left\| \left\| \frac{f \left( \varphi_x^{2^{-n}} \right)}{2^{- n \alpha}} \right\|_{L^p \left( x \in \mathbb{R}^d \right)} \right\|_{\ell^q \left( n \right)} < + \infty .
			\end{equation}
			
			The same statement holds for the local Besov spaces $\mathcal{B}_{p, q, \mathrm{loc}}^{\alpha} ( \mathbb{R}^d )$, when one replaces $L^p ( x \in \mathbb{R}^d )$ by $L^p \left( x \in K \right)$ for all compact $K \subset \mathbb{R}^d$ in the condition above.
		
		\item\label{item:prop_equivalent_besov_norm_2} If $\alpha > 0$ and $\alpha \notin \mathbb{N}$, then $f \in \mathcal{B}_{p, q}^{\alpha} ( \mathbb{R}^d )$ if and only if for all $0 \leq \left| k \right| < \alpha$, $\partial^k f \in L^p$ and for any $h_0 > 0$,
			\begin{equation}
				\left\| \left\| \frac{\partial^k f \left( x + h \right) - \sum\limits_{0 \leq \left| l \right| < \alpha - \left| k \right|} \frac{1}{l !} \partial^{k + l} f \left( x \right) h^k}{h^{\alpha - \left| k \right|}} \right\|_{L^{p} \left( x \right)} \right\|_{L_h^q \left( h \in B \left( 0, h_0 \right) \right)} < + \infty .
				\label{eq:Besov_Taylor_bound}
			\end{equation}
			
			The same statement holds for the local Besov spaces $\mathcal{B}_{p, q, \mathrm{loc}}^{\alpha} ( \mathbb{R}^d )$, when one replaces $L^p ( x \in \mathbb{R}^d )$ by $L^p \left( x \in K \right)$ for all compact $K \subset \mathbb{R}^d$ in the condition above.
	\end{enumerate}
\end{proposition}

We choose to provide a proof of this proposition for the sake of completeness, although we believe that these properties are well-known in the literature of Besov spaces. 
For instance, \autoref{item:prop_equivalent_besov_norm_1} is proven to be equivalent to the usual \enquote{Littlewood-Paley} definition of Besov spaces in \cite[Corollary~1.12]{MR2250142}.
Also, see \cite{MR500920, MR4099476} for examples of papers using a definition of Besov spaces similar to \autoref{item:prop_equivalent_besov_norm_2}.

The techniques used in the proof below are very reminiscent of those used in the remainder of this paper, which is another reason for us to include it.

\begin{proof}
We prove the assertions separately.
The local versions of the results are established with similar calculations so we only prove the global versions.
	\begin{enumerate}
		\item[\ref{item:prop_equivalent_besov_norm_1}] The direct implication is straightforward. Now let us concentrate on the converse. Let $f \in \mathcal{D}^{\prime} ( \mathbb{R}^d )$ and $\varphi \in \mathcal{D} ( \mathbb{R}^d )$ be a test-function as in the statement, we shall show that $f \in \mathcal{B}_{p, q}^{\alpha}$. We \enquote{tweak} the test-function $\varphi$ as in \autoref{LEMMA:TWEAK}. Let $r \in \mathbb{N}$ be such that $r > - \alpha$, and fix distinct $\lambda_0, \cdots \lambda_{r-1}$ small enough so that we can define $\hat{\varphi} \in \mathscr{B}_{\left\lfloor \alpha \right\rfloor}^{r}$ as in \autoref{LEMMA:TWEAK}. Note that also:
			\begin{equation}
				\left\| \left\| \frac{f \left( \hat{\varphi}_x^{2^{-n}} \right)}{2^{- n \alpha}} \right\|_{L^p \left( x \right)} \right\|_{\ell^q \left( n \right)} < + \infty .
			\end{equation}
			
			As above, set $\check{\varphi}:=\hat{\varphi}^{1/2}-\hat{\varphi}^{2}$ so that by mollification we have the following decomposition for all $\psi \in \mathcal{D} ( \mathbb{R}^d )$ (see \autoref{DEF:APPROX}):
				\begin{align}
					f \left( \psi_x^{2^{-n}} \right) & = \int_{\mathbb{R}^d} f \left( \hat{\varphi}_z^{2^{- n}} \right) \left( \hat{\varphi}^{2^{- n + 1}} * \psi^{2^{-n}} \right) \left( z - x \right) d z \\
					& \quad + \sum\limits_{m \geq n} \int_{\mathbb{R}^d} f \left( \hat{\varphi}_z^{2^{- m}} \right) \left( \check{\varphi}^{2^{- m}} * \psi^{2^{- n}} \right) \left( z - x \right) d z .
				\end{align}
				
			Substituting $\tilde{z} \coloneqq z - x$ and integrating only on the support of the integrand:
				\begin{align}
					f \left( \psi_x^{2^{-n}} \right) & = \int_{\tilde{z} \in B \left( 0, 2^{- n + 1} \right)} f \left( \hat{\varphi}_{\tilde{z} + x}^{2^{- n}} \right) \left( \hat{\varphi}^{2^{- n + 1}} * \psi^{2^{-n}} \right) \left( \tilde{z} \right) d \tilde{z} \\
					& \quad + \sum\limits_{m \geq n} \int_{\tilde{z} \in B \left( 0, 2^{- n + 1} \right)} f \left( \hat{\varphi}_{\tilde{z} + x}^{2^{- m}} \right) \left( \check{\varphi}^{2^{- m}} * \psi^{2^{- n}} \right) \left( \tilde{z} \right) d \tilde{z} .
				\end{align}
				
			Now we use the estimates (see \autoref{LEMMA:CONV_INEQ}):
				\begin{equation}
					\begin{dcases}
						\left\| \hat{\varphi}^{2^{- n + 1}} * \psi^{2^{-n}} \right\|_{L^{\infty}} & \leq 2^{n d} \left\| \psi \right\|_{L^{\infty}} \left\| \hat{\varphi} \right\|_{L^1} \leq 2^{n d} \left\| \psi \right\|_{\mathcal{C}^r} \left\| \hat{\varphi} \right\|_{L^1}  , \\
						 \left\| \check{\varphi}^{2^{- m}} * \psi^{2^{- n}} \right\|_{L^{\infty}} & \leq 2^{n \left( r + d \right) - m r} \left\| \psi \right\|_{\mathcal{C}^r} \left\| \check{\varphi} \right\|_{L^1} . \\
					\end{dcases}
				\end{equation}
			
			This yields:
				\begin{align}
					\sup\limits_{\psi \in \mathscr{B}_{}^{r}} \left| f \left( \psi_x^{2^{-n}} \right) \right| & \leq 2^{n d} \left\| \hat{\varphi} \right\|_{L^1} \int_{\tilde{z} \in B \left( 0, 2^{- n + 1} \right)} f \left( \hat{\varphi}_{\tilde{z} + x}^{2^{- n}} \right) d \tilde{z} \\
					& \quad + \sum\limits_{m \geq n} 2^{n \left( r + d \right) - m r} \left\| \check{\varphi} \right\|_{L^1}  \int_{\tilde{z} \in B \left( 0, 2^{- n + 1} \right)} f \left( \hat{\varphi}_{\tilde{z} + x}^{2^{- m}} \right) d \tilde{z} . 
				\end{align}
			
			In order to simplify notations, denote $B_n \coloneqq B \left( 0, 2^{- n} \right)$ Thus, integrating over $x$:
				\begin{align}
					\left\| \sup\limits_{\psi \in \mathscr{B}_{}^{r}} \left| f \left( \psi_x^{2^{-n}} \right) \right| \right\|_{L^p \left( x \right)} & \leq 2^{n d} \left\| \hat{\varphi} \right\|_{L^1} \int_{\tilde{z} \in B_{n - 1}} \left\| f \left( \hat{\varphi}_{\tilde{z} + x}^{2^{- n}} \right) \right\|_{L^p \left( x \right)} d \tilde{z} \\
					& \quad + \sum\limits_{m \geq n} 2^{n \left( r + d \right) - m r} \left\| \check{\varphi} \right\|_{L^1}  \int_{\tilde{z} \in B_{n - 1}} \left\| f \left( \hat{\varphi}_{\tilde{z} + x}^{2^{- m}} \right) \right\|_{L^p \left( x \right)} d \tilde{z} . 
				\end{align}
			
			In those integrals in $\tilde{z}$, the integrand is actually constant so that after integration, we obtain:
				\begin{align}
					\left\| \sup\limits_{\psi \in \mathscr{B}_{}^{r}} \left| f \left( \psi_x^{2^{-n}} \right) \right| \right\|_{L^p \left( x \right)} & \lesssim \left\| f \left( \hat{\varphi}_{x}^{2^{- n}} \right) \right\|_{L^p \left( x \right)} + \sum\limits_{m \geq n} 2^{\left( n - m \right) r} \left\| f \left( \hat{\varphi}_{x}^{2^{- m}} \right) \right\|_{L^p \left( x \right)} .
				\end{align}
			
			Then:
				\begin{align}
					\left\| f \right\|_{\mathcal{B}_{p, q}^{\alpha}} & \lesssim \left\| \left\| \frac{f \left( \hat{\varphi}_{x}^{2^{- n}} \right)}{2^{- n \alpha}} \right\|_{L^p \left( x \right)} \right\|_{\ell^q \left( n \right)} + \left\| \sum\limits_{m \geq n} 2^{\left( n - m \right) \left( r + \alpha \right)} \left\| \frac{f \left( \hat{\varphi}_{x}^{2^{- m}} \right)}{2^{- m \alpha}} \right\|_{L^p \left( x \right)} \right\|_{\ell^q \left( n \right)} .
				\end{align}
		Since we chose $r + \alpha > 0$, applying Jensen's inequality then interverting the sums in $m$ and $n$ in the second term of the right-hand side yields as announced $\left\| f \right\|_{\mathcal{B}_{p, q}^{\alpha}} < + \infty$.
		
		\item[\ref{item:prop_equivalent_besov_norm_2}] For simplicity, we reason with $h_0 = 1$ (but the same arguments generalise to any $h_0 > 0$).
		Let us first concentrate on the direct statement. Let $f \in \mathcal{B}_{p, q}^{\alpha}$. By definition of the Besov space and the distributional definition of $\partial^k f$, it holds that for $0 \leq  \left| k \right| < \alpha$, $\partial^k f \in \mathcal{B}_{p, q}^{\alpha - \left| k \right|} ( \mathbb{R}^d )$, so it suffices to prove the claim for $\left| k \right| = 0$.
			First, let us show that $f \in L^p$, in the sense that there exists $\tilde{f} \in L^p$ such that $f = \tilde{f}$ as distributions.
			We reason by mollification.
			Fix $\varphi \in \mathscr{B}^{r}$ a single test-function such that $\int \varphi = 1$ and $\int x^l \varphi \left( x \right) d x = 0$ for $1 \leq \left| l \right| < \alpha$. For $m, n \in \mathbb{N}$, define $\tilde{f}_{m, n} \left( x \right) \coloneqq f ( \varphi_x^{2^{- m - n}} )$.
			Using the embedding $\ell^q \subset \ell^{\infty}$, it holds that:
				\begin{align}
					\left\| \left\| \frac{\tilde{f}_{m, n} \left( x \right) - \tilde{f}_{m, n + 1} \left( x \right)}{2^{- n \alpha}} \right\|_{L^p \left( x \right)} \right\|_{l^{\infty} \left( n \right)} & \leq \left\| \left\| \frac{\tilde{f}_{m, n} \left( x \right) - \tilde{f}_{m, n + 1} \left( x \right)}{2^{- n \alpha}} \right\|_{L^p \left( x \right)} \right\|_{l^{q} \left( n \right)} \\
					& = \left\| \left\| \frac{f \left( \varphi_x^{2^{- m - n}} - \varphi_x^{2^{- m - n - 1}} \right)}{2^{- n \alpha}} \right\|_{L^p \left( x \right)} \right\|_{l^{q} \left( n \right)} .
				\end{align}	
					
			Let $\psi \coloneqq \varphi - \varphi^{\frac{1}{2}}$. Note that $ \frac{1}{C}\psi\in  \mathscr{B}_{\left\lfloor \alpha \right\rfloor}^{r}$ for a suitable $C>0$. In particular, we have:
				\begin{align}
					\left\| \left\| \frac{f \left( \varphi_x^{2^{- m - n}} - \varphi_x^{2^{- m - n - 1}} \right)}{2^{- n \alpha}} \right\|_{L^p \left( x \right)} \right\|_{l^{q} \left( n \right)} & = \left\| \left\| \frac{f \left( \psi_x^{2^{- m - n}} \right)}{2^{- n \alpha}} \right\|_{L^p \left( x \right)} \right\|_{l^{q} \left( n \right)} \\
					& \leq 2^{- m \alpha} \left\| f \right\|_{\mathcal{B}_{p, q}^{\alpha}} . 
				\end{align}
			We deduce that $\| \tilde{f}_{m, n} \left( x \right) - \tilde{f}_{m, n + 1} \left( x \right) \|_{L^p \left( x \right)} \leq 2^{- \left( m + n \right) \alpha} \| f \|_{\mathcal{B}_{p, q}^{\alpha}}$.
			This implies that for each $m \in \mathbb{N}$, the sequence $( \tilde{f}_{m, n})_{n \in \mathbb{N}}$ is Cauchy in $L^p$. Hence it has a limit, which we call $\tilde{f}_m \in L^p$, satisfying:
				\begin{equation}
					\left\| \tilde{f}_m - \tilde{f}_{m, n} \right\|_{L^p} \xrightarrow[n \to + \infty]{} 0 .
				\end{equation}
				
			By summation of a geometric series, we even have $\| \tilde{f}_m - \tilde{f}_{m, n} \|_{L^p} \lesssim 2^{- \left( m + n \right) \alpha}$.
			And for any $n_0 \in \mathbb{N}$, the following series converges in $L^p$:	
				\begin{equation}
					\tilde{f}_m = \tilde{f}_{m, n_0} + \sum\limits_{n = n_0}^{+ \infty} \left( \tilde{f}_{m, n + 1} - \tilde{f}_{m, n} \right) .
				\end{equation}
			
			Also, since for any $m, n \in \mathbb{N}$, $\tilde{f}_{m, n} = \tilde{f}_{m + 1, n - 1}$, we deduce by triangle inequality that:
				\begin{align}
					\left\| \tilde{f}_{m + 1} - \tilde{f}_m \right\|_{L^p} & \leq \left\| \tilde{f}_{m + 1} - \tilde{f}_{m + 1, n - 1} \right\|_{L^p} + \left\| \tilde{f}_{m} - \tilde{f}_{m, n} \right\|_{L^p} \\
					& \xrightarrow[n \to + \infty]{} 0 .
				\end{align}
			Thus for all $m \in \mathbb{N}$, $\tilde{f}_{m} = \tilde{f}_{m + 1} \eqqcolon \tilde{f}$, where this equality holds in $L^p$  (hence also in $\mathcal{D}^{\prime}$) .
			Now let us show that $\tilde{f} = f$ in $\mathcal{D}^{\prime}$.
			Let $\psi \in \mathcal{D} ( \mathbb{R}^d )$. By mollification, 
				\begin{equation}
					\left( \tilde{f} - f \right) \left( \psi \right) = \lim\limits_{n \to + \infty} \int_{\mathbb{R}^d} \left( \tilde{f} \left( x \right) - f \left( \varphi_x^{2^{-n}} \right) \right) \psi \left( x \right) d x .
				\end{equation}
			But for $n \in \mathbb{N}$, H\"older's inequality:
				\begin{align}
					\left| \int_{\mathbb{R}^d} \left( \tilde{f} \left( x \right) - f \left( \varphi_x^{2^{-n}} \right) \right) \psi \left( x \right) d x \right| & \lesssim \left\| \tilde{f} \left( x \right) - f \left( \varphi_x^{2^{-n}} \right) \right\|_{L^p \left( x \right)} \\
					& = \left\| \tilde{f} - \tilde{f}_{0, n} \right\|_{L^p} \\
					& \xrightarrow[n \to + \infty]{} 0 .
				\end{align}
			This establishes the announced equality.
			Now we establish \eqref{eq:Besov_Taylor_bound} for $\left| k \right| = 0$, which, according to our previous remark, suffices to establish \eqref{eq:Besov_Taylor_bound}. For $x, h \in \mathbb{R}^d$, set $T_f^{\alpha} \left( x, h \right) \coloneqq f \left( x + h \right) - \sum_{0 \leq \left| l \right| < \alpha} \frac{1}{l !} \partial^{l} f \left( x \right) h^l$, so that we shall show:
			\begin{equation}
		\left\| \left\| \frac{T_f^{\alpha} \left( x, h \right)}{\left| h \right|^{\alpha}} \right\|_{L^{p} \left( x \right)} \right\|_{L_h^q \left( h \right)} < + \infty .
			\end{equation}
			
			For each $h \in \mathbb{R}^d$ we consider $m_h \in \mathbb{Z}$ defined to be explicited later. 
			We write with the notations of the previous item:
				\begin{align}
					T_f^{\alpha} \left( x, h \right) & = T_{f_{m_h}}^{\alpha} \left( x, h \right) \\ 
					& = T_{f_{m_h, 0}}^{\alpha} \left( x, h \right) + \left( T_{f_{m_h}}^{\alpha} \left( x, h \right) - T_{f_{m_h, 0}}^{\alpha} \left( x, h \right) \right) .
				\end{align}
			
			More explicitly:
				\begin{align}
					T_f^{\alpha} \left( x, h \right) & = f \left( \varphi_{ x + h }^{2^{- m_h}} - \sum\limits_{0 \leq \left| l \right| < \alpha} \frac{\left( -1 \right)^{\left| l \right|} h^l}{l !} \partial^l \left( \varphi_x^{2^{- m_h}} \right) \right) \\
					& \quad + \sum\limits_{n = 0}^{+ \infty} f \left( \varphi_{x + h}^{2^{- m_h - n - 1}} - \varphi_{x + h}^{2^{- m_h - n}} \right) \\
					& \quad - \sum\limits_{0 \leq \left| l \right| < \alpha} \sum\limits_{n = 0}^{+ \infty} \frac{h^l}{l !} \partial^l f \left( \varphi_{x}^{2^{- m_h - n - 1}} - \varphi_{x}^{2^{- m_h - n}} \right) .
				\end{align}
			
			Now we bound each of these terms using our definition of Besov spaces. For $h \in B \left( 0, 1 \right)$ and $z \in \mathbb{R}^d$, define:
				\begin{equation}
					\begin{dcases}
						\psi \left( z \right) & \coloneqq \varphi \left( z - 2^{m_h} h \right) - \sum\limits_{0 \leq \left| l \right| < \alpha} \frac{\left( - 2^{m_h} h \right)^l}{l !} \partial^l \varphi \left( z \right) , \\
						\tilde{\psi} \left( z \right) & \coloneqq \varphi^{\frac{1}{2}} \left( z \right) - \varphi \left( z \right) .
					\end{dcases}
				\end{equation}
			Then:
				\begin{align}
					\left\| \left\| \frac{T_f^{\alpha} \left( x, h \right)}{\left| h \right|^{\alpha}} \right\|_{L^{p} \left( x \right)} \right\|_{L_h^q \left( h \right)} & \leq \left\| \left\| \frac{f \left( \psi_x^{2^{- m_h}} \right)}{\left| h \right|^{\alpha}} \right\|_{L^{p} \left( x \right)} \right\|_{L_h^q \left( h \right)} \\
					& \quad + \sum\limits_{n = 0}^{+ \infty} \left\| \left\| \frac{f \left( \tilde{\psi}_{x + h}^{2^{- m_h - n}} \right)}{\left| h \right|^{\alpha}} \right\|_{L^{p} \left( x \right)} \right\|_{L_h^q \left( h \right)} \\
					& \quad + \sum\limits_{0 \leq \left| l \right| < \alpha} \frac{1}{l !} \sum\limits_{n = 0}^{+ \infty} \left\| \left\| \frac{ \partial^l f \left( \tilde{\psi}_{x}^{2^{- m_h - n}}\right)}{\left| h \right|^{\alpha - \left| l \right|}} \right\|_{L^{p} \left( x \right)} \right\|_{L_h^q \left( h \right)} .
				\end{align}
				
			Changing variable $\tilde{x} = x + h$ in the second term and noting that there exists a constant $C > 0$ such that $\frac{1}{C}\psi, \frac{1}{C}\tilde{\psi} \in  \mathscr{B}_{\left\lfloor \alpha \right\rfloor}^{r}$ (note that actually $\supp \left( \psi \right) \subset B \left( 0, 2 \right)$ rather than $B \left( 0, 1 \right)$, but this is not a problem after invoking the definition of Besov spaces \autoref{def:Besov_spaces_by_local_means} for $n_0 = -1$; note also that this is where we require $\alpha$ to be non-integer), we obtain:
				\begin{align}
					& \left\| \left\| \frac{T_f^{\alpha} \left( x, h \right)}{\left| h \right|^{\alpha}} \right\|_{L^{p} \left( x \right)} \right\|_{L_h^q \left( h \right)} \\
					& \quad \leq C \left\| \left\| \sup\limits_{\psi \in \mathscr{B}_{\left\lfloor \alpha \right\rfloor}^{r}} \left| \frac{f \left( \psi_x^{2^{- m_h}} \right)}{\left| h \right|^{\alpha}} \right\|_{L^{p} \left( x \right)} \right| \right\|_{L_h^q \left( h \right)} \\
					& \quad \quad + \sum\limits_{0 \leq \left| l \right| < \alpha} \frac{2 C}{l !} \sum\limits_{n = 0}^{+ \infty} \left\| \left\| \sup\limits_{\tilde{\psi} \in \mathscr{B}_{\left\lfloor \alpha \right\rfloor}^{r}} \left| \frac{ \partial^l f \left( \tilde{\psi}_{x}^{2^{- m_h - n}}\right)}{\left| h \right|^{\alpha - \left| l \right|}} \right| \right\|_{L^{p} \left( x \right)} \right\|_{L_h^q \left( h \right)} .
				\end{align}
			
			To conclude, it is enough to prove that if $f \in \mathcal{B}_{p, q}^{\alpha}$, then:
				\begin{equation}
					v \left( f \right) \coloneqq \left\| \left\| \sup\limits_{\psi \in \mathscr{B}_{\left\lfloor \alpha \right\rfloor}^{r}} \left| \frac{f \left( \psi_x^{2^{- m_h}} \right)}{\left| h \right|^{\alpha}} \right| \right\|_{L^{p} \left( x \right)} \right\|_{L_h^q \left( h \right)} < + \infty .
				\end{equation}
			
			We cut the integral in $h$ along the annuli: for $n \in \mathbb{N}$, set $B_n \coloneqq B \left( 2^{- \left( n + 1 \right)}, 2^{- n} \right)$ then:
				\begin{equation}
					v \left( f \right) = \left\| \left\| \left\| \sup\limits_{\psi \in \mathscr{B}_{\left\lfloor \alpha \right\rfloor}^{r}} \left| \frac{f \left( \psi_x^{2^{- m_h}} \right)}{\left| h \right|^{\alpha}} \right| \right\|_{L^{p} \left( x \right)} \right\|_{L_h^q \left( h \in B_n \right)} \right\|_{\ell^q \left( n \right)} .
				\end{equation}
				
			Now we choose $m_h$ so that for $h \in B_n, m_h = n$.
			Using the fact that $\left\| \frac{1}{\left| h \right|^{\alpha}} \right\|_{L_h^q \left( h \in B_n \right)} \lesssim 2^{n \alpha}$ uniformly in $n \in \mathbb{N}$, we get $v \left( f \right) \lesssim \left\| f \right\|_{\mathcal{B}_{p, q}^{\alpha}} < + \infty$, which concludes the direct statement of the proposition.
		
		Now let us turn to the converse. Assume that 	for all $0 \leq \left| k \right| < \alpha$, $\partial^k f \in L^p$ and that \eqref{eq:Besov_Taylor_bound} holds, we shall prove that $f \in \mathcal{B}_{p, q}^{\alpha}$. On the one hand, when $\psi \in \mathscr{B}^{r}$ it holds that $f \left( \psi_x \right) = \int_{\supp \left( \psi \right)} f \left( y - x \right) \psi \left( y \right) d y$, so that:
			\begin{align}
				\left\| \sup\limits_{\psi \in \mathscr{B}^{r}} \left| f \left( \psi_x \right) \right| \right\|_{L^p \left( x \right)} & \leq \left\| \sup\limits_{\psi \in \mathscr{B}^{r}} \left\| \psi \right\|_{\mathscr{B}^{r}} \int_{B \left( 0, 1 \right)} \left| f \left( y - x \right) \right| d y \right\|_{L^p \left( x \right)} \\
				& \leq \int_{B \left( 0, 1 \right)} \left\| f \right\|_{L^p} d y \\
				& < + \infty .
			\end{align}
		
		On the other hand, when $\psi \in \mathscr{B}_{\left\lfloor \alpha \right\rfloor}^{r}$, by subtracting a suitable Taylor polynomial it holds that $f \left( \psi_x^{2^{-n}} \right) = \int_{\mathbb{R}^d} T_{f}^{\alpha} \left( x, y - x \right) \psi \left( \frac{y - x}{2^{-n}} \right) 2^{nd} d y$, so that:
			\begin{align}
				& \left\| \left\| \sup\limits_{\psi \in \mathscr{B}_{\left\lfloor \alpha \right\rfloor}^{r}} \left| \frac{f \left( \psi_x^{2^{-n}} \right)}{2^{- n \alpha}} \right| \right\|_{L^p \left( x \right)} \right\|_{\ell^q \left( n \right)} \\
				& \quad \leq \left\| \left\| \sup\limits_{\psi \in \mathscr{B}_{\left\lfloor \alpha \right\rfloor}^{r}} \left| \frac{\int_{B \left( x, 2^{-n} \right)} \left| T_{f}^{\alpha} \left( x, y - x \right) \right| \left\| \psi \right\|_{\mathscr{B}_{\left\lfloor \alpha \right\rfloor}^{r}} d y}{2^{- n \left( \alpha + d \right)}} \right| \right\|_{L^p \left( x \right)} \right\|_{\ell^q \left( n \right)} \\
				& \quad \leq \left\| \int_{B \left( 0, 2^{-n} \right)} \left( \frac{\left\| T_{f}^{\alpha} \left( x, h \right) \right\|_{L^p \left( x \right)}}{\left| h \right|^{\alpha}} \right) \left( \frac{\left| h \right|^{\alpha}}{2^{- n \left( \alpha + d \right)}} \right) d h \right\|_{\ell^q \left( n \right)} .
			\end{align}
		Since $\sup\limits_{n \in \mathbb{N}} \int_{B \left( 0, 2^{-n} \right)} \frac{\left| h \right|^{\alpha}}{2^{- n \left( \alpha + d \right)}} d h < + \infty$, we have by applying Jensen's inequality as well as switching summations, 
			\begin{align}
				& \left\| \int_{B \left( 0, 2^{-n} \right)} \left( \frac{\left\| T_{f}^{\alpha} \left( x, h \right)  \right\|_{L^p \left( x \right)}}{\left| h \right|^{\alpha}} \right) \left( \frac{\left| h \right|^{\alpha}}{2^{- n \left( \alpha + d \right)}} \right) d h \right\|_{\ell^q \left( n \right)} \\
				& \quad \lesssim \left( \sum\limits_{n = 0}^{+ \infty} \int_{B \left( 0, 2^{-n} \right)} \left( \frac{\left\| T_{f}^{\alpha} \left( x, h \right) \right\|_{L^p \left( x \right)}}{\left| h \right|^{\alpha}} \right)^q \left( \frac{\left| h \right|^{\alpha}}{2^{- n \left( \alpha + d \right)}} \right) d h \right)^{\frac{1}{q}} \\
				& \quad = \left( \int_{B \left( 0, 1 \right)} \sum\limits_{n = 0}^{\left\lfloor - \log_2 \left( \left| h \right| \right) \right\rfloor} \left( \frac{\left\| T_{f}^{\alpha} \left( x, h \right)  \right\|_{L^p \left( x \right)}}{\left| h \right|^{\alpha}} \right)^q \left( \frac{\left| h \right|}{2^{- n}} \right)^{\alpha + d} \frac{d h}{\left| h \right|^{d}} \right)^{\frac{1}{q}} .
			\end{align}
		Since $\sup\limits_{h \in B \left( 0, 1 \right)} \sum\limits_{n = 0}^{\left\lfloor - \log_2 \left( \left| h \right| \right) \right\rfloor} \left( \frac{\left| h \right|}{2^{- n}} \right)^{\alpha + d} < + \infty$, we deduce:
		\begin{align}
				\left\| \left\| \sup\limits_{\psi \in \mathscr{B}_{\left\lfloor \alpha \right\rfloor}^{r}} \left| \frac{f \left( \psi_x^{2^{-n}} \right)}{2^{- n \alpha}} \right| \right\|_{L^p \left( x \right)} \right\|_{\ell^q \left( n \right)} & \lesssim \left\| \left\| \frac{T_{f}^{\alpha} \left( x, h \right)}{\left| h \right|^{\alpha}} \right\|_{L^{p} \left( x \right)} \right\|_{L_h^q \left( h \right)} \\
				& < + \infty .
			\end{align}
		Hence, $f \in \mathcal{B}_{p, q}^{\alpha}$.
	\end{enumerate}
\end{proof}

\section{A technical lemma on series}

In this section, we establish the following technical result, used in \autoref{section:proof_of_corollary} for the proof of \autoref{corollary:reconstruction_with_coherence_homogeneity}, and in \autoref{section:application_young_multiplication} for the proof of \autoref{prop:Young_multiplication_besov}.

\begin{lemma}\label{LEMMA:tech}
	Let $(f_k:\R^d\rightarrow \R_+)_{k\in \N}$ be a family of positive functions, and $\left( a_{k, n} \right)_{k, n \in \mathbb{N}} \in \mathbb{R}_+^{\mathbb{N}^2}$ be a sequence of positive reals.
	Consider the sequence defined by:
	\begin{equation*}
		u_n \coloneqq \sum\limits_{k = 0}^{+ \infty} a_{k, n} \int_{|x|\leq  2^{-k + 1} } 2^{k d} f_k \left( x \right) dx \quad \text{for $n \in \mathbb{N}$} . 
	\end{equation*} 
	Assume that there exists $A > 0$ such that:
			\begin{equation}\label{EQN:tech_lemma_bnds}
				\begin{dcases}
					\text{for all } n \in \mathbb{N} :& \sum\limits_{k = 0}^{+ \infty} a_{k, n} \leq A , \\
					\text{for all } k \in \mathbb{N} :& \sum\limits_{n = 0}^{+ \infty} a_{k, n} \leq A .
				\end{dcases}
			\end{equation}
	
	 Fix $q\in \N$ and assume also that either of the following conditions is satisfied:
		\begin{enumerate}[(i), ref =\emph{(\roman*)}]
			\item\label{item:technical_lemma_condition_1} $\left\| \left\| f_k \left( x \right) \right\|_{\ell^\infty \left( k \right)} \right\|_{L_x^q \left( x \in B \left( 0, 2 \right) \right)} < + \infty$, or
			\item\label{item:technical_lemma_condition_2} There exists $\mu \in \ell^{q_1} \left( n \in \mathbb{N} \right)$, $\nu \in L_h^{q_2} \left( h \in B ( 0, 2 ) \right)$ for some $q_1, q_2 \in \left[ 1, + \infty \right]$ satisfying $\frac{1}{q} = \frac{1}{q_1} + \frac{1}{q_2}$, and such that $f_k \left( x \right) \leq \mu \left( k \right) \nu \left( x \right)$.
		\end{enumerate}		
	Then $(u_n)_{n\in \N} \in \ell^q$.
	\label{lemma:technical_result}
\end{lemma}

\begin{proof}
Let us first prove the result under the assumption \ref{item:technical_lemma_condition_1}.
	For simplicity, we assume that $q < + \infty$ since the case for $q=\infty$ is straightforward.
	For a fixed $n \in \mathbb{N}$, we apply Jensen's inequality on $u_n$ and we obtain:
	\begin{align*}
	\left (\frac{\left| u_n \right|}{\sum_{k=0}^{\infty}a_{k,n}}\right )^q &\leq \frac{1}{\sum_{k=0}^{\infty}a_{k,n}}\sum\limits_{k = 0}^{+ \infty} a_{k, n} \left( \int_{|x|\leq 2^{-k + 1}} 2^{k d} f_k \left( x \right) dx \right)^q.
	\end{align*}
	Hence, we have
	\begin{align}\label{EQN:u_nq_2}
	| u_n|^q &\leq \left (\sum\limits_{k = 0}^{+ \infty} a_{k, n}\right )^{q-1} \sum\limits_{k = 0}^{+ \infty} a_{k, n}\left( \int_{|x|\leq 2^{-k + 1}} 2^{k d} f_k \left( x \right) dx \right)^q \\
	&\leq A^{q-1} \sum\limits_{k = 0}^{+ \infty} a_{k, n}\left( \int_{|x|\leq 2^{-k + 1}} 2^{k d} f_k \left( x \right) dx \right)^q
	\end{align}
	from \eqref{EQN:tech_lemma_bnds}. 

Applying Jensen's inequality, on the integral, with the probability measure:
	\begin{displaymath}
		c_{d}2^{kd} \mathds{1}_{|x|\leq  2^{-k + 1}} dx,	\quad \text{where $c_d= \frac{1}{\text{Vol}(B(0,2))}$}, 
	\end{displaymath}
		we obtain
	\begin{displaymath}
	\left (\int_{|x|\leq 2^{-k + 1}} 2^{k d} f_k \left( x \right) dx\right )^q \leq c_d^{1-q}\int_{|x|\leq 2^{-k + 1}} 2^{k d} \left|f_k \left(x\right)\right|^q dx .
	\end{displaymath}
	which yields:
	\begin{displaymath}
	\|(u_n)_{n\in \N}\|_{\ell^q}^q \leq c_d^{1-q}A^{q-1} \sum_{n=0}^{\infty}\sum\limits_{k = 0}^{+ \infty} a_{k, n} \int_{|x|\leq 2^{-k + 1} } 2^{k d} |f_k \left( x \right)|^q dx .
	\end{displaymath}
	
	Note that the integral can be decomposed over annuli:
	\begin{displaymath}
		\int_{|x|\leq 2^{-k + 1} } 2^{k d} |f_k \left( x \right)|^q dx=\sum_{l=k - 1}^{\infty}\int_{2^{-(l+1)}\leq |x|\leq 2^{-l} }  2^{k d} |f_k \left( x \right)|^q dx.
	\end{displaymath}
	By applying Tonelli's theorem and rearranging sequences we obtain:
	\begin{align*}
	\left\| u \right\|_{\ell^q}^q & \leq c_d^{1-q}A^{q-1} \sum\limits_{n = 0}^{+ \infty} \sum\limits_{k = 0}^{+ \infty} \sum\limits_{l = k - 1}^{+ \infty} a_{k, n} \int_{2^{-(l+1)}\leq |x|\leq 2^{-l} } 2^{k d} |f_k \left( x \right)|^q dx \\
	& = c_d^{1-q}A^{q-1} \sum\limits_{n = 0}^{+ \infty} \sum\limits_{0 \leq k \leq l + 1 < + \infty} a_{k, n} \int_{_{2^{-(l+1)}\leq |x|\leq 2^{-l} } } 2^{k d} |f_k \left( x \right)|^q dx \\\nonumber
	& = c_d^{1-q}A^{q-1} \sum\limits_{l = -1}^{+ \infty} \sum\limits_{k = 0}^{l + 1} \sum\limits_{n = 0}^{+ \infty} a_{k, n} \int_{_{2^{-(l+1)}\leq |x|\leq 2^{-l} } } 2^{k d} |f_k \left( x \right)|^q dx \\
	& = c_d^{1-q}A^{q-1} \sum\limits_{l = -1}^{+ \infty} \int_{2^{-(l+1)}\leq |x|\leq 2^{-l} } \sum\limits_{k = 0}^{l + 1} 2^{k d} |f_k \left( x \right)|^q \sum\limits_{n = 0}^{+ \infty} a_{k, n} dx .
	\end{align*}
	
	Applying assumption \eqref{EQN:tech_lemma_bnds} we obtain \begin{displaymath}
	\left\| u \right\|_{\ell^q}^q \leq c_d^{1-q}A^{q} \sum\limits_{l = -1}^{+ \infty} \int_{2^{-(l+1)}\leq |x|\leq 2^{-l} } \sum\limits_{k = 0}^{l + 1} 2^{k d} |f_k \left( x \right)|^q dx .
	\end{displaymath}
	
	For $l$ fixed, we can obtain the following
	\begin{displaymath}
		\sum\limits_{k = 0}^{l+1} 2^{k d} |f_k \left( x \right)|^q \leq \frac{2^{dl+2d}-1}{2^d-1}\|f_{k}(x)\|_{l^\infty}^q\leq \frac{2^{2d}}{2^d -1} 2^{ld}\|(f_{k}(x))_{k\in\N}\|_{l^\infty}^q.
	\end{displaymath}
	 Hence:
	\begin{align*}
	\left\| u \right\|_{\ell^q}^q & \leq c_{d}^{1-q} A^q \frac{2^{2d}}{2^d -1} \sum\limits_{l = -1}^{+ \infty} \int_{2^{-(l+1)}\leq |x|\leq 2^{-l} } 2^{l d} \left\| f_k \left( x \right) \right\|_{\ell^\infty \left( k \right)}^{q} dx \\
	& \leq  c_{d}^{1-q}A^q \frac{2^{2d}}{2^d -1}\sum\limits_{l = -1}^{+ \infty} \int_{2^{-(l+1)}\leq |x|\leq 2^{-l} } \left\| f_k \left( x \right) \right\|_{\ell^\infty \left( k \right)}^{q} \frac{d x}{\left| x \right|^d}\\
	&=c_{d}^{1-q}A^q \frac{2^{2d}}{2^d -1} \int_{0\leq |x|\leq 2} \left\| f_k \left( x \right) \right\|_{\ell^{\infty} \left( k \right)}^{q} \frac{d x}{\left| x \right|^d}\\
	&=c_{d}^{1-q}A^q \frac{2^{2d}}{2^d -1} \left\| \left\| f_k \left( x \right) \right\|_{\ell^\infty \left( k \right)} \right\|_{L_x^q \left( x \in B \left( 0, 2 \right) \right)}^q \\
	& < + \infty ,
	\end{align*}
	from \ref{item:technical_lemma_condition_1}.
	
	Now let us prove the result under the assumption \ref{item:technical_lemma_condition_2}.
	By Jensen's inequality:
	\begin{equation}
		\left| u_n \right|^q \lesssim \sum\limits_{k = 0}^{+ \infty} a_{k, n} \left| \mu \left( k \right) \right|^q \int_{|x|\leq  2^{-k + 1} } 2^{k d} \left| \nu \left( x \right) \right|^q dx .
	\end{equation}

We sum over $n \in \mathbb{N}$ and intervert summations in $k$ and $n$ by Fubini:
	\begin{equation}
		\sum\limits_{n \in \mathbb{N}} \left| u_n \right|^q \lesssim \sum\limits_{k = 0}^{+ \infty} \left( \sum\limits_{n = 0}^{+ \infty} a_{k, n} \right) \left| \mu \left( k \right) \right|^q \int_{|x|\leq  2^{-k + 1} } 2^{k d} \left| \nu \left( x \right) \right|^q dx .
	\end{equation}

By assumption, the sum of $a_{k, n}$ is bounded, hence:
	\begin{equation}
		\sum\limits_{n \in \mathbb{N}} \left| u_n \right|^q \lesssim \sum\limits_{k = 0}^{+ \infty} \left| \mu \left( k \right) \right|^q \int_{|x|\leq  2^{-k + 1} } 2^{k d} \left| \nu \left( x \right) \right|^q dx .
	\end{equation}

We now apply H\"older's inequality with the conjugate exponents $\frac{q_1}{q}$, $\frac{q_2}{q}$:
	\begin{equation}
		\sum\limits_{n \in \mathbb{N}} \left| u_n \right|^q \lesssim \left( \sum\limits_{k = 0}^{+ \infty} \left| \mu \left( k \right) \right|^{q \frac{q_1}{q}} \right)^{\frac{q}{q_1}} \left( \sum\limits_{k = 0}^{+ \infty} \left( \int_{|x|\leq  2^{-k + 1} } 2^{k d} \left| \nu \left( x \right) \right|^q dx \right)^{\frac{q_2}{q}} \right)^{\frac{q}{q_2}} .
	\end{equation}

Applying Jensen's in the integral:
	\begin{equation}
		\left\| u \right\|_{\ell^q}^q \lesssim \left\| \mu \right\|_{\ell^{q_1}}^q \left( \sum\limits_{k = 0}^{+ \infty}  \int_{|x|\leq  2^{-k + 1} } 2^{k d} \left| \nu \left( x \right) \right|^{q_2} dx \right)^{\frac{q}{q_2}} .
	\end{equation}
	
We decompose the domain of the integral as an union of dyadic annuli:
	\begin{equation}
		\left\| u \right\|_{\ell^q}^q \lesssim \left\| \mu \right\|_{\ell^{q_1}}^q \left( \sum\limits_{k = 0}^{+ \infty} \sum\limits_{l = k - 1}^{+ \infty} \int_{ 2^{- \left( l + 1 \right)} \leq |x|\leq  2^{-l} } 2^{k d} \left| \nu \left( x \right) \right|^{q_2} dx \right)^{\frac{q}{q_2}} .
	\end{equation}
	
Interverting the sums:	
	\begin{equation}
		\left\| u \right\|_{\ell^q}^q \lesssim \left\| \mu \right\|_{\ell^{q_1}}^q \left( \sum\limits_{l = -1}^{+ \infty} \int_{ 2^{- \left( l + 1 \right)} \leq |x|\leq  2^{-l} } \left( \sum\limits_{k = 0}^{l + 1} 2^{k d} \right) \left| \nu \left( x \right) \right|^{q_2} dx \right)^{\frac{q}{q_2}} .
	\end{equation}
	
Now in this integral $\sum\limits_{k = 0}^{l + 1} 2^{k d} \lesssim \frac{1}{\left| x \right|^d}$ and thus:
	\begin{equation}
		\left\| u \right\|_{\ell^q}^q \lesssim \left\| \mu \right\|_{\ell^{q_1}}^q \left\| \nu \right\|_{L_x^{q_2} \left( x \in B ( 0, 2 ) \right)}^{q} .
	\end{equation}

By assumption, this is finite, and thus our assertion is proved.
\end{proof}

\subsection*{Acknowledgements}
We are very grateful to Lorenzo Zambotti for suggesting the topic of this investigation and for his helpful discussions and insights. 
We also acknowledge Cyril Labb\'e and Francesco Caravenna for their thorough review of our earlier drafts.
We thank Florian Bechtold, Peter Friz, David Pr\"omel and Federico Sclavi for their useful feedback.  

This project has received funding from the
European Union's Horizon 2020 research and innovation programme under the
Marie Sk\l{}odowska-Curie grant agreement No 754362.

\begin{figure}[ht]
	\centering
	\includegraphics[height=10mm]{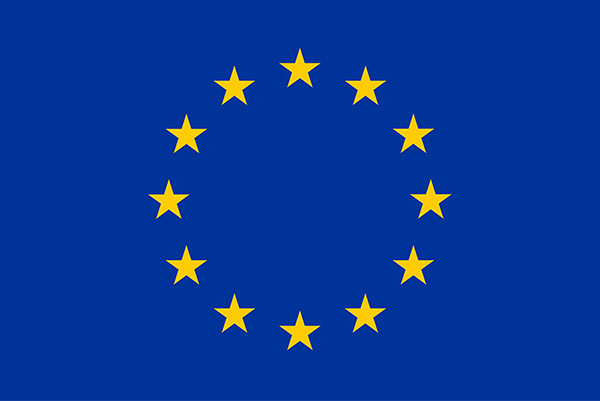}
\end{figure}

%
%
%
%

\bibliographystyle{siam}
\bibliography{biblio}

\end{document}